\documentclass[preprint]{elsarticle}

\usepackage{lineno,hyperref}
\modulolinenumbers[5]

\usepackage{graphicx,amssymb,amsmath,amsthm,color,amsmath}
\usepackage{color,multirow}
\usepackage{mathrsfs}
\usepackage{setspace}
\usepackage{adjustbox}
\usepackage{array}
\usepackage{algorithmic,setspace}
\usepackage[ruled]{algorithm}
\usepackage{longtable}
\usepackage{caption}

\def\ind{\mathbf{ind}}
\def\col{\mathbf{col}}
\def\tcol{\text{col}}
\def\h{\mathbf{h}}
\def\g{\mathbf{g}}
\def\j{\mathbf{j}}
\def\myphi{\mathbf{phi}}
\def\mypsi{\mathbf{psi}}
\def\tphi{\text{phi}}
\def\tpsi{\text{psi}}

\def\cq{\{q\}}
\def\vf{\mathbf{f}}
\def\vc{\mathbf{c}}
\def\cl{\boldsymbol{\ell}}
\def\vfa{\vf^{\rm{a}}}
\def\vfc{\vf^{\rm{c}}}
\def\vfr{\mathbf{f}^{\rm{r}}}
\def\vd{\mathbf{d}}
\def\vD{\mathbf{D}}
\def\vDw{\mathbf{D}^W}
\def\vDc{\mathbf{D}^C}
\def\Mc{M_c}
\def\Mw{M_w}

\def\SRB{\mathrm{SR}_{\mathrm{B}}}
\def\SRD{\mathrm{SR}_{\mathrm{D}}}

\def\fontalg{} 

\def\Nq{N_b}
\def\kq{k(q)}

\def\be{\begin{equation}}
\def\ee{\end{equation}}
\def\R{\mathbb{R}}
\def\SR{\mathrm{SR}}
\def\sr{\mathrm{sr}}
\def\name{\mathrm{namef}}
\def\pars{\mathrm{pars}}
                         
\def\PRD{\mathrm{PRD}}

\def\tol{\mathrm{tol}}
\def\prd{{\mathrm{prd}}}
\def\prdo{{\mathrm{prd}_0}}

\DeclareCaptionFormat{algor}{ 
  \hrulefill\par\offinterlineskip\vskip1pt 
    \textbf{#1#2}#3\offinterlineskip\hrulefill}
\DeclareCaptionStyle{algori}{singlelinecheck=off,format=algor,labelsep=space}
\newcommand{\SWITCH}[1]{\STATE \textbf{switch} (#1)}
\newcommand{\ENDSWITCH}{\STATE \textbf{end switch}}
\newcommand{\CASE}[1]{\STATE \textbf{case} #1\textbf{:} \begin{ALC@g}}
\newcommand{\ENDCASE}{\end{ALC@g}}
\newcommand{\OTHERWISE}{\STATE \textbf{otherwise} \begin{ALC@g}}
\newcommand{\ENDOTHERWISE}{\end{ALC@g}}

\newdefinition{rem}{Remark}
\newdefinition{ex}{Example}

\begin{document}

\begin{frontmatter}

\title{ Construction of wavelet dictionaries for ECG modelling  }

\author{Dana \v{C}ern{\'a} } 
\address{ Department of Mathematics and Didactics of Mathematics, Technical University of Liberec, Studentsk\'a 2, 461 17 Liberec, Czech Republic }

\author{Laura Rebollo-Neira}
\address{ Mathematics Department, Aston University, B4 7ET, Birmingham, UK}

\begin{abstract}

{\em{Background and Objective:}} The purpose of sparse modelling of ECG signals is to represent an ECG  record, given by sample points, as a linear combination of as few elementary components as possible. This can be achieved by creating a redundant set, called a  dictionary, from where the elementary components are selected.  The success in  
sparsely representing an ECG record  depends on the nature of the dictionary being considered.  In this paper we focus on the construction of different  families of wavelet dictionaries, which  are appropriate  for the purpose of reducing dimensionality of  ECG signals through sparse representation modelling. 

{\em{Method:}} The suitability of wavelet dictionaries for ECG modelling, applying the 
Optimized Orthogonal Matching Pursuit  approach  for the selection process, 
was demonstrated in a previous work on the MIT-BIH Arrhythmia database consisting
of 48 records  each of which of 30 min length. This paper complements the previous one by presenting the technical details, methods, algorithms, and MATLAB software facilitating the
construction of different families of wavelet dictionaries. The implementation allows for straightforward further extensions to include additional wavelet families.

{\em{Results:}} The sparsity in the representation of an ECG record significantly improves
in relation to the sparsity produced by the corresponding wavelet basis. This result holds true for the 17 wavelet families considered here.

{\em{Conclusions:}} Wavelet dictionaries contribute to  the representation of an 
ECG record as a superposition of fewer components than those needed by the 
wavelet basis. The software for the construction of wavelet dictionaries, 
which has been made available to support the material in this paper, could be of
assistance to a broad range of application relying on dimensionality reduction as
a first step of further ECG signal analysis.

\end{abstract}


\end{frontmatter}


\section{Introduction}
\label{introduction}

The electrocardiogram (ECG) is a routine 
test for clinical medicine. It plays a crucial role in 
the diagnosis of a broad range of anomalies in the human 
heart; from arrythmias to myocardial infarction. 

The widely available digital ECG data has facilitated 
the development of algorithms for  
 ECG processing and interpretation. 
In particular, the literature for computerized
 arrhythmia detection and classification is   
 extensive. Useful review matterial \cite{LSC16,LMM18} 
 can help with the 
 introduction to state of the art 
 techniques, which nonetheless keeps growing 
 \cite{AFA17, HRH19, ANP18, LYM17, Pla18}. 
 
A common first step in ECG modeling consists in  
reducing the dimensionality of the signal. This 
entails to represent the informational content of the 
record by means of significantly fewer parameters than 
the number of samples in the digital ECG.
When the aim is to reproduce the original signal 
at low level distortion, the 
step is frequently realized through transformations 
such as the Wavelet Transform 
and the Discrete Cosine Transform. 
In the last few years alternative approaches, 
falling within the category 
of sparse representation of ECG signals,
have been considered. 

Within the sparse representation framework, an ECG record
 is represented as a linear combination of elementary 
components, called 
{\em{atoms,}}  which are selected from a redundant set, 
called a {\em{dictionary}}. The success of the 
methods developed within this framework depends on 
both, the selection technique and the proposed dictionary. 
The selection techniques which are widely applied 
for sparse representation of general signals 
are either greedy pursuit strategies 
\cite{MZ92,PRK93,RNL02}, or strategies based on 
 minimization of the 1-norm
as a cost function \cite{CDS98}. Suitable dictionaries 
 depend on the class of signals 
being processed. These can be designed at 
hoc or be learned from training data.   
Sparse representation of ECG signals has been 
tackled by both these approaches, e.g. 
\cite{AGL15} learns dictionaries using some part of the 
records for ECG compression and 
\cite{RR18} uses Gabor dictionaries for 
structuring features for classification.

In a recent publication \cite{RNC19} we have 
shown that wavelet dictionaries, derived  from 
known wavelet families, are suitable for  representing 
an ECG record as a linear combination of fewer elementary
 components than those required by a wavelet basis. 
The model was shown to be successful for
dimensionality reduction and lossy compression.
As far as compression is concerned the method advanced 
in \cite{RNC19} produces compression results
improving upon previously reported benchmarks
\cite{LKL11, MZD15, PSS16, TZW18} for the MIT-BIH 
 Arrhythmia data set without pre-processing. 
 With regard to dimensionality reduction, 
 wavelet dictionaries considerably improve upon the
 results achieved with the wavelet basis of the
 same family \cite{RNC19,LRN19}. This result motivated the 
 present Communication. While in \cite{RNC19} the 
dictionaries have been used to demonstrate their suitability 
for dimensionality reduction of ECG signals at low 
level distortion, the details of their  
 numerical construction were not given.  
This paper complements the previous work by   
presenting the algorithms for building dictionaries 
 from the following mother wavelet prototypes: 
\begin{enumerate}
    \item[1)] Chui-Wang linear spline wavelet \cite{CW92} 
    \item[2)] Chui-Wang quadratic spline wavelet \cite{CW92}
    \item[3)] Chui-Wang cubic spline wavelet \cite{CW92}
    \item[4)] primal CDF97 wavelet \cite{CDF92}
    \item[5)] dual CDF97 wavelet \cite{CDF92}
    \item[6)] primal CDF53 wavelet \cite{CDF92}
    \item[7)] linear spline wavelet with short support and 2 vanishing moments \cite{Chen95, Han06}
    \item[8)] quadratic spline wavelet with short support and 3 vanishing moments \cite{Chen95, Han06}
    \item[9)] cubic spline wavelet with short support and 4 vanishing moments \cite{Chen95, Han06}
    \item[10)] Daubechies wavelet with $3$ vanishing moments \cite{Daub88}
    \item[11)] Daubechies wavelet with $4$ vanishing moments \cite{Daub88}
    \item[12)] Daubechies wavelet with $5$ vanishing moments \cite{Daub88}
    \item[13)] symlet with $3$ vanishing moments \cite{Daub93}
    \item[14)] symlet with $4$ vanishing moments \cite{Daub93}
    \item[15)] symlet with $5$ vanishing moments \cite{Daub93}
    \item[16)] coiflet with 2 vanishing moments and support of length 6 which is the most regular \cite{Daub93}
    \item[17)] coiflet with 3 vanishing moments and support of length 8 \cite{Daub93}
\end{enumerate} 
 
The method proposed in \cite{RNC19} for modelling 
a given ECG signal proceeds as follows. 
 Assuming that the signal is given as an $N$-dimensional 
 array, this array is partitioned into $Q$ cells 
 $\vfc\cq,q=1,\ldots,Q$. 
Thus, each cell $\vfc\cq$ is an $\Nq$-dimensional vector, 
which is modeled by an atomic decomposition of the form
\be
\label{atom}
\vfa\cq = \sum_{n=1}^{\kq}\vc\cq(n) \vd_{\ell\cq(n)}.
\ee
For each cell $q$, the atoms $\vd_{\ell_{\cq(n)}},\,n=1,\ldots,\kq$ are selected from a dictionary through the greedy Optimized Orthogonal Matching Pursuit  (OOMP) algorithm \cite{RNL02,LRN16}.
The array $\cl\cq$ is a vector whose components  
 $\ell\cq(n),\,n=1,\ldots,\kq$ 
contain the indices of the selected atoms 
for decomposing the $q$-th cell in the signal 
 partition. 
The OOMP method, for selecting these indices and 
 computing the corresponding coefficients $\vc\cq ,\,n=1,\ldots,\kq$ in 
\eqref{atom}, is fully implemented by the OOMP function 
included as a tool of the software.  
 
Each of the proposed dictionaries consists 
of two components.   
One of the components contains  a few  
  elements, say $\Mc$, from a discrete cosine basis. 
This component of the dictionary allows for the fact that 
ECG signals are normally superimposed to a 
smooth background. It is given as a $\Nq \times \Mc$ 
 matrix $\vDc$.
The other component is the wavelet-based dictionary, 
which is given as a $\Nq \times \Mw$ matrix $\vDw$. 
Thus, the whole dictionary $\vD$ is an $\Nq \times (\Mc+\Mw)$ matrix obtained by the horizontal concatenation of 
$\vDc$ and $\vDw$.  The next section is dedicated to 
the construction of $\vDw$.

The paper is organized as follows. 
Sec. 2 gives  all the details for the construction of different wavelet prototypes and the concomitant wavelet dictionaries generated by those prototypes. Secs 3 and 4 deliver details and examples demonstrating the use of the MATLAB software for modelling ECG signals within the proposed framework.  

The software has been made available on a dedicated webpage \cite{RNC19a}.
The implementation allows for straightforward further extension of the options for wavelet types.

\section{Method}
\label{sec_wavelet_dict}

In this section  we produce  all the pseudo-codes for the construction of  wavelets dictionaries, which can  be  used  to achieve the model of every segment in a signal partition.  As already mentioned,   each dictionary is obtained by taking the prototypes from a  wavelet basis and translating them  within a shorter step  than  that corresponding to the wavelet basis.
    
Throughout the paper we adopt the following notation. Boldface fonts are used to indicate Euclidean vectors and matrices.  Standard mathematical fonts  are used to indicate components,  e.g., $\vd \in \R^N$ is a vector of $N$-components $d(i) \in \R,\, i=1,\ldots,N$ and $\vD \in \R^{N\times M}$  is a matrix of  elements  $D(i,j),\, i= 1\ldots, N, j=1,\ldots, M$. The symbol $L^2(\R)$ denotes the space of square integrable functions.

Wavelets are usually constructed starting from a {\it multiresolution analysis}, which is 
a sequence $\left\{ V_j \right\}_{j = j_0}^{\infty}$ of closed subspaces of the space $L^2\left(\mathbb{R}\right)$ which are nested and their union is dense in $L^2\left(\mathbb{R}\right)$, 
i.e.,
\begin{equation}
V_j \subset V_{j+1} \quad \forall \, j \geq j_0, \quad \overline{\bigcup_{j = j_0}^{\infty} V_j }=L^2\left(\mathbb{R}\right).
\end{equation}
We assume that there exists a function $\phi \in L^2 \left( \mathbb{R} \right)$ such that for $j \geq j_0$ functions
\begin{equation} \label{scaling_functions}
\phi_{j,k} \left( x \right)  = 2^{j/2} \phi  \left( 2^j x - k \right), \quad k \in \mathbb{Z}, 
\end{equation}
form uniformly stable bases of the spaces $V_j$, i.e., the bases are Riesz bases with bounds independent of the level $j$, see e.g. \cite{Cohen2003}. The functions $\phi_{j,k}$ are called {\it scaling functions} and the function $\phi$ is called a {\it generator} of scaling functions. Next we present a method for the actual construction of the scaling functions. 

\subsection{Generation of scaling functions}

We assume that $\phi$ has a compact support $\left[ 0, K \right]$ for some $K \in \mathbb{N}$. From the nestedness of the multiresolution spaces $V_j$, it follows that there exists a {\it scaling filter} $\mathbf{h} =  \left( h(1), \ldots, h(K+1) \right)$ such that 
\begin{equation} \label{scaling_equation}
\phi\left(x\right) = \sum_{k=1}^{K+1} h \left( k \right) \phi \left( 2x + 1 -k\right) \quad \forall \ x \in \mathbb{R}.
\end{equation}
If $\int_0^K \phi \left( x \right) dx = c \neq 0$
then, integrating (\ref{scaling_equation}), we obtain
\begin{equation}
    c = \sum_{k=1}^{K+1} h \left( k \right) \frac{c}{2}
\end{equation}
which implies that $\mathbf{h}$ has to be normalized such that
\begin{equation} \label{normalization_h}
    \sum_{k=1}^{K+1} h \left( k \right) = 2.
\end{equation}

The scaling equation~(\ref{scaling_equation}) enables computing values of the scaling function $\phi$ at points $k/2^u$ for $k=0, \ldots, K 2^u$, $u \in \mathbb{N}$. First we compute values of $\phi$ at integer points.
Since $\text{supp} \, \phi = \left[ 0, K \right]$, we have $\phi \left( k \right)=0$ for $k \neq \left( 0, K \right)$. Let us define a vector
\begin{equation} \label{def_Phi}
\Phi=\left( \phi \left( 0 \right), \ldots, \phi \left( K-1 \right) \right)^T,
\end{equation}
where the $(.)^T$ indicates the
transpose operation.
Substituting $x=0, \ldots, K-1,$ into (\ref{scaling_equation}), we obtain 
\begin{eqnarray} \label{phi_at_integers}
\Phi \left( i \right) &=& \phi \left( i - 1 \right) =
\sum_{k=1}^{K+1} h \left( k \right) \phi \left( 2i - 1 -k\right) \\
\nonumber &=& \sum_{j=2i-1}^{2i-K-1} h \left( 2i -j  \right) \phi \left( j -1 \right) = \sum_{j=2i-1}^{2i-K-1} h \left( 2i -j  \right) \Phi \left( j \right).
\end{eqnarray}
We set $h(k)=0$ for $k<1$ and $k>K+1$ and define a matrix $\mathbf{A}$ by
\begin{equation} \label{def_matrix_A} 
    A \left( i, j \right)= h \left( 2i -j \right), \quad i,j=1, \ldots, K.
\end{equation}
Then, (\ref{phi_at_integers}) is equivalent to
\begin{equation}
    \Phi = \mathbf{A} \Phi.
\end{equation}
This means that $\Phi$ is an eigenvector corresponding to the eigenvalue $1$ of the matrix $\mathbf{A}$. If the multiplicity of this eigenvalue
is $1$, then $\Phi$ is given uniquely up to a multiplication by a constant.
Our aim is to compute a vector $\myphi$ such that
\begin{equation} \label{def_tphi}
\tphi \left( m \right) = \phi \left( \frac{m-1}{2^u} \right),
\quad m=1, \ldots, K \, 2^u +1,
\end{equation}
for a chosen level $u \in \mathbb{N}$.
From (\ref{def_Phi}) and (\ref{def_tphi}) we have
\begin{equation} \label{set_phi_integers}
    \tphi \left( 2^u \left( l - 1 \right) + 1 \right) = \phi \left( l-1 \right) = \Phi \left( l \right), \quad l=1, \ldots, K.
\end{equation}
We compute values of $\phi$ at points $l/2$. Note that for $l$ even we already know these values.
Using (\ref{scaling_equation}) and (\ref{set_phi_integers}) we obtain
\begin{eqnarray}
\tphi \left( l 2^{u-1} + 1\right) &=& \phi\left( \frac{l}{2} \right) = \sum_{k=1}^{K+1} h \left( k \right) \phi \left( l + 1 -k\right)   \\
\nonumber &=& \sum_{k=1}^{K+1} h \left( k \right) \tphi \left( \left( l + 1 -k\right) 2^u +1 \right) 
\end{eqnarray}
for $l=1, 3, \ldots, 2 N -1$. Similarly, we compute values of $\phi$ at points $l/4$, and thus we continue until we determine values at points $l/2^u$. More precisely, for $i=1, \ldots ,u$ we assume that we know values of $\phi$ at $l/ 2^{i-1}$, $l=0, \ldots, K 2^{i-1} + 1,$ and we compute the values 
\begin{equation} \label{m_to_x}
\tphi \left( m \right) = \phi \left( x \right), \quad  m=x 2^u + 1, \quad 
x=l/ 2^{i-1}+1/ 2^i. 
\end{equation}
Using (\ref{scaling_equation}) we obtain
\begin{eqnarray} \label{phi_m}
\tphi \left( m \right) &=& \phi\left( x \right) = \sum_{k=1}^{K+1} h \left( k \right) \phi \left( 2 x + 1 -k\right)   \\
\nonumber &=& \sum_{k=1}^{K+1} h \left( k \right) \tphi \left( \left( 2x + 1 - k\right) 2^u +1 \right).
\end{eqnarray}

\begin{rem} Some scaling functions such as spline scaling functions are known in an explicit form and their values can be evaluated directly. However, an advantage of our approach is that it is more general and can be used for a large class of wavelet families.
\end{rem}

\subsection{Construction of wavelet generators from scaling functions} 

Let $W_j$ be complement spaces such that $V_{j+1}= V_j \oplus W_j$, where $\oplus$ denotes a direct sum. 
Wavelet functions $\psi_{j,k}$ are constructed in the form:
\begin{equation} 
\psi_{j,k} \left( x \right) = 2^{j/2} \psi \left( 2^j x - k \right),
\quad k \in \mathbb{Z},
\end{equation}
to be a basis for $W_j$ and such that  
\begin{equation} \label{def_wavelet_basis}
\mathcal{B} =\left\{ \phi_{j_0,k}, k \in \mathbb{Z} \right\} \cup \left\{ \psi_{j,k}, k \in \mathbb{Z}, j \geq j_0 \right\}
\end{equation}
called a {\it wavelet basis}, is a Riesz basis of the space $L^2 \left( \mathbb{R} \right)$. 

Since $W_j \subset V_{j+1}$ there exists a vector $\mathbf{g}=\left( g(1), \ldots, g(M+1) \right)$
such that
\begin{equation} \label{wavelet_equation}
 \psi \left( x \right) = \sum_{k=1}^{M+1} g \left( k \right) \phi \left( 2x + 1 -k \right).  
\end{equation}
The vector $\mathbf{g}$ is called a {\it wavelet filter}. From (\ref{wavelet_equation}) we have
\begin{equation}
\text{supp} \, \psi = \left[ 0, \frac{M+K}{2}  \right].
\end{equation}
In Algorithm~\ref{algor_WaveletGen} we compute a vector $\mypsi$ such that 
\begin{equation} \label{def_tpsi}
\tpsi \left( m \right) = \psi \left( \frac{m-1}{2^u} \right), \quad m=1, \ldots, \left( M+K \right) 2^{u-1} +1, 
\end{equation} in the following way. Due to (\ref{wavelet_equation}) and (\ref{def_tpsi}), we have
\begin{eqnarray} 
\nonumber \tpsi \left( m \right) &=& \psi \left( \frac{m-1}{2^u} \right) = 
\sum_{k=1}^{M+1} g \left( k \right) \phi \left( \left(m-1 \right) \, 2^{1-u} + 1 -k \right) \\
& = & \sum_{k=1}^{M+1} g \left( k \right)  \tphi \left( 2 m -1 + \left( 1-k \right) \, 2^u \right).
\end{eqnarray}
The sum in the last equation is computed as a cyclic sum.
For 
\begin{equation}
m=1, \ldots, ( M + K ) 2^{u-1} + 1
\end{equation}
we set $\tpsi \left( m \right)=0$ and for $k=1, \ldots, M+1$ we do
\begin{equation}
\tpsi \left( m \right) = \tpsi  \left( m \right) + g \left( k \right)  \tphi \left( 2 m -1 + \left( 1-k \right) \, 2^u \right), 
\end{equation}
if $1 \leq 2 m -1 + \left( 1-k \right) \, 2^u \leq K 2^u + 1$. 
Using the substitution 
\begin{equation}
2 m -1 + \left( 1-k \right) \, 2^u \ = 2i -1,
\end{equation}
for $i=1, \ldots,  K 2^{u-1} + 1$, we obtain
\begin{equation} \label{compute_psi}
\tpsi \left( i + \left( k - 1 \right) 2^{u-1} \right) = \tpsi  \left(  i + \left( k - 1 \right) 2^{u-1} \right) + g \left( k \right)  \tphi \left( 2 i -1 \right). 
\end{equation}

Algorithm~\ref{algor_WaveletGen} computes vectors $\myphi$ and $\mypsi$ for given scaling and wavelet filters.
The  filters corresponding to  the wavelet families supported by the software are
given in Appendix A (Algorithm~\ref{algor_Filters}).

\begin{center}
{\fontalg 
  \captionsetup{type=figure}
  \captionof{algorithm}{ \\
  Procedure 
[$\myphi$,$\mypsi$] = WaveletGen($\h$,$\g$,u) } \label{algor_WaveletGen}

\begin{algorithmic}

\STATE{ {\bf Input:} }\, 

\begin{tabular}{ p{1cm} p{9cm}  }
$\h$ & scaling filter   \\
$\g$ & wavelet filter   \\
$u$ & level (integer) that determines points $l/2^u$ \\
\end{tabular}

\STATE{ {\bf Output:}}\, 

\begin{tabular}{ p{1cm} p{9cm}  }
$\myphi$ & vector of scaling function values (c.f. (\ref{def_tphi})) \\
$\mypsi$ & vector of wavelet values (c.f. (\ref{def_tpsi})) \\
\end{tabular}

\vspace{2mm}

\STATE{ $K=\text{length} (\h) -1$}
\COMMENT{support length of $\phi$}

\STATE{ $\h = 2 \, \h /  \text{sum} \left( \h \right) $ } \COMMENT{normalization of $\mathbf{h}$ (c.f. (\ref{normalization_h}))}

\STATE{\COMMENT{Compute a matrix $\mathbf{A}$ using (\ref{def_matrix_A})}}
\STATE{ $\mathbf{A} = \text{zeros} (K)$}

\FOR { $i=1:K$ }

\FOR { $j=1:K$ }

\IF {$ 1 \leq 2i-j \leq K+1$}

\STATE{$A \left( i, j \right) = h \left( 2i-j \right)$}

\ENDIF

\ENDFOR

\ENDFOR

\STATE{\COMMENT{Compute eigenvalues and eigenvectors of the matrix $\mathbf{A}$} }

\STATE{[V,D]=eig($\mathbf{A}$)} 

\STATE{\COMMENT{Find an index of a column corresponding to eigenvalue 1} }

\STATE{$k=0$} \COMMENT{$k$ is the multiplicity of eigenvalue 1}  
\FOR{ $i=1:K$ }

\IF{ $ \left| (D(i,i)-1 \right| <10^{-7} $} 
       
\STATE{ column=$i$; $k=k+1$ }

\ENDIF 

\ENDFOR

\IF{ $ k \neq 1$ }
   
\STATE{ error(`Impossible to construct scaling function: eigenvalue 1 must have multiplicity 1') }

\ELSE
   
\STATE{ $\myphi= \text{zeros}(K 2^u+1, 1) $ }

\STATE{\COMMENT{Eigenvector $V(:,\text{column})$ represents values of $\phi$ at integer points} }

\STATE{ $\tphi( 1:2^u:(K-1) 2^u+1) = V(:,\text{column})$ } 
\COMMENT{c.f. (\ref{set_phi_integers})}

\COMMENT{Compute values of $\phi$ at points $l/2^u$}

\FOR{ $i = 1 : u$ } 

\FOR{ $l=1:K  2^{i-1}$ } 

\STATE{ $x=2^{-i}+(l-1) \, 2^{-i+1}$;  $m=x \, 2^u+1$}
\COMMENT{c.f. (\ref{m_to_x})}
           
\FOR{ $k=1:K+1$ }
               
               \IF{ $ 0 \leq 2x-k+1 \leq K$}
               
                  \STATE{ $\tphi(m)=\tphi(m)+ h(k) \tphi((2x-k+1) \, 2^u+1)$ }
                 \COMMENT{c.f. (\ref{phi_m})}
                  
               \ENDIF 
              
\ENDFOR

\ENDFOR 

\ENDFOR

\STATE{$M= \text{length} (\g)-1$}

\COMMENT{Compute $\mypsi$ containing values of $\psi$ at points $l/2^u$}

\STATE{$ \mypsi = \text{zeros}((K+M) \, 2^{u-1}+1,1) $ }

\FOR{ $k=1:M+1$ }  

\STATE{ $i_1=(k-1) \, 2^{u-1}+1$;  $i_2=(k-1) \, 2^{u-1}+1+K \, 2^{u-1}$ }
    
\STATE{ $\tpsi(i_1:i_2)= \tpsi(i_1:i_2)+g(k) \, \tphi(1:2:K \, 2^u+1)$ } 
\COMMENT{c.f. (\ref{compute_psi})}

\ENDFOR

\ENDIF

\hrulefill 
\end{algorithmic}
}
\end{center}

\subsection{Construction of  wavelet bases  and dictionaries} 

Hereafter we drop all normalization factors and normalize all the vectors once they have been constructed.
Note that in (\ref{def_wavelet_basis}) we used a translation parameter $k \in \mathbb{Z}$ and since $\mathcal{B}$ is a Riesz basis the functions from $\mathcal{B}$ are linearly independent. 
Now, we choose a parameter $b$ such that $b = 2^{-m}$ for some integer $m$. We define functions
\begin{equation} \label{scaling_functions_dictionary}
\phi_{j_0,k,b} \left( x \right)  =  \phi  \left( 2^{j_0} x - b k \right), \quad k \in \mathbb{Z}, 
\end{equation}
and
\begin{equation} \label{wavelets_dictionary}
\psi_{j,k,b} \left( x \right) =  \psi \left( 2^j x - b k \right), \quad
k \in \mathbb{Z}, \quad j \geq j_0,
\end{equation}
which form a redundant dictionary \cite{ARN05, ARN08, RNX10}. Obviously, $b=1$ corresponds to a basis.

\begin{figure} [!ht] 
 	\includegraphics[height=.2\textheight]{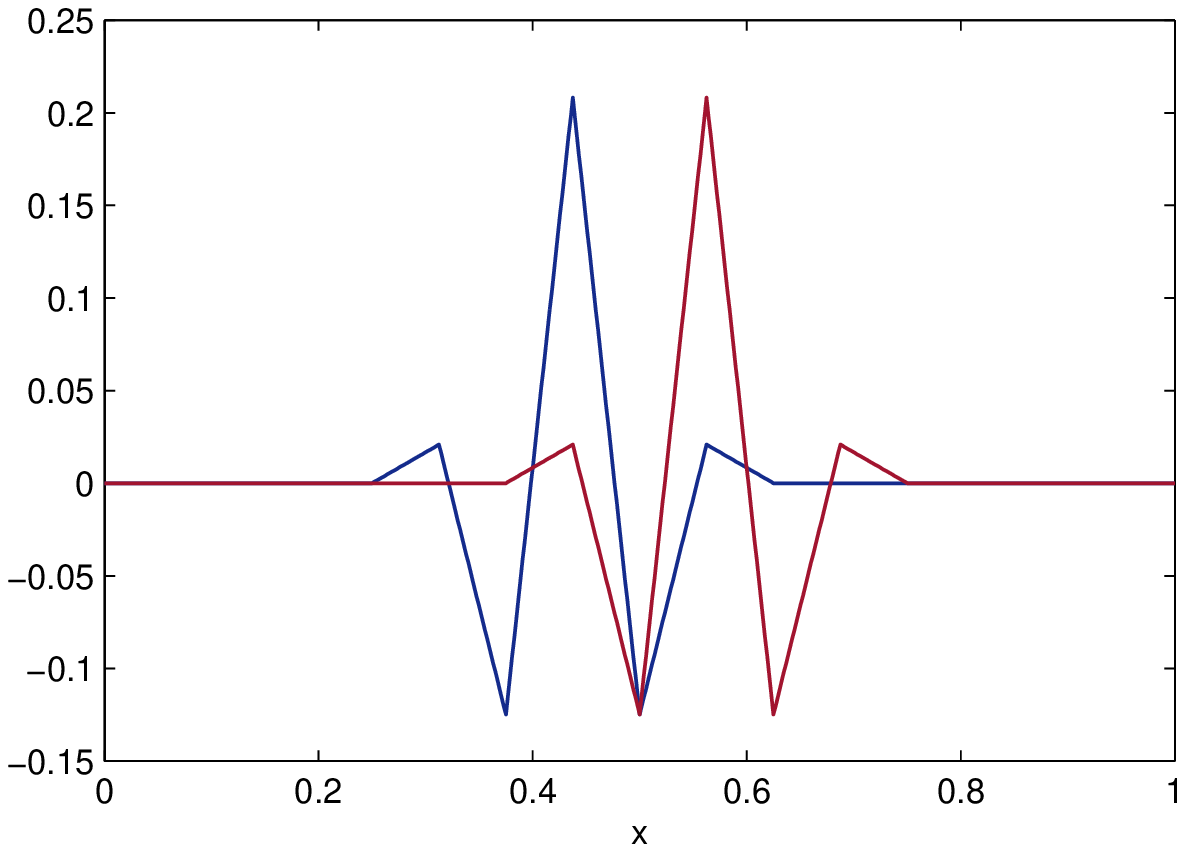}
 	\includegraphics[height=.2\textheight]{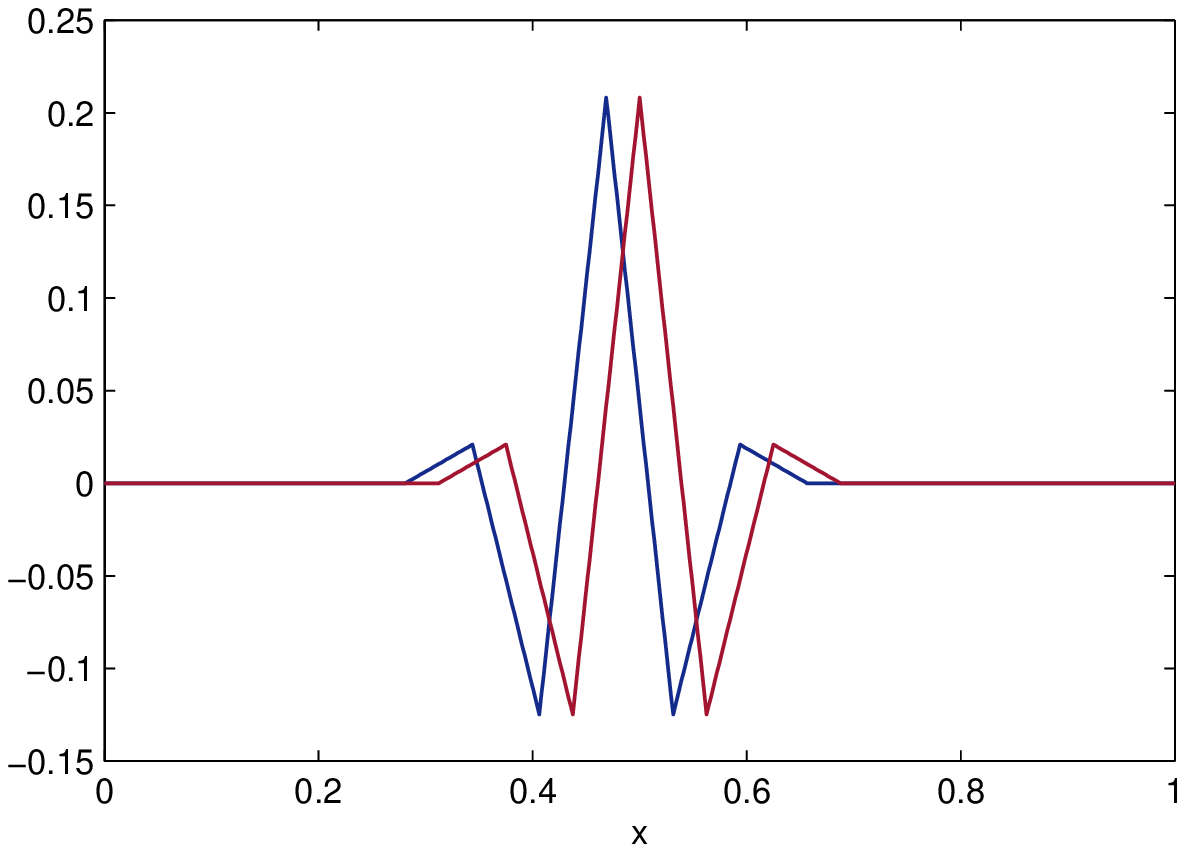}
  \caption{ Wavelet functions taken from a basis (left) and a dictionary (right) corresponding to a linear spline-wavelet prototype from \cite{CW92}.}
  \label{figure1}
\end{figure}

The left graph of Figure~\ref{figure1} shows two consecutive
wavelet functions taken from a linear spline bases \cite{CW92}.  The right  graph of Figure~\ref{figure1} 
corresponds to two  consecutive  wavelet functions taken from the dictionary spanning the same space which  corresponds to $b=1/4$.

Algorithm~\ref{algor_WaveletDict} constructs a discrete dictionary, i.e., a matrix  $\vDw$ which contains values of functions from (\ref{scaling_functions_dictionary}) and (\ref{wavelets_dictionary}) at $N_b$ equidistant points for some chosen levels determined by the vector $\mathbf{j}$. Since Algorithm~\ref{algor_WaveletGen} enables us to construct values at points of the form $l/2^u$, we evaluate functions
(\ref{scaling_functions_dictionary}) and (\ref{wavelets_dictionary}) at the points
\begin{equation} \label{Nb_points}
x=\frac{l}{2^r}, \quad l=0, \ldots, N_b-1, \quad r= \left\lceil \frac{\log \left( N_b -1 \right)}{\log \left( 2 \right)} \right\rceil,
\end{equation}
where $\left\lceil y \right\rceil$ denotes the smallest integer number larger than $y$.

For a chosen vector of levels $\j$, we define a vector of indices $\ind$ such that $ \text{ind}(1)$ is the number of scaling functions at level $j(1)$, and $\text{ind}(l)$ is the number of wavelets at level $j(l-1)$ for $l=1, \ldots, J$, where $J$ is the length of $\j$.
We have
\begin{equation}
    \text{supp} \, \phi_{j(1),k,b} = \left[ \frac{bk}{ 2^{j(1)}}, \frac{bk+K}{ 2^{j(1)} } \right],
    \quad \text{supp} \, \psi_{j,k,b} = \left[ \frac{bk}{ 2^{j} }, \frac{ bk+ \frac{K+M}{2} }{ 2^j } \right].
\end{equation}
Comparing the supports of these functions and the interval 
\begin{equation} \label{def_I}
I=\left[ 0, \frac{ N_b -1 }{2^r } \right]
\end{equation}
which contains the points from (\ref{Nb_points}), we find that the number of inner scaling functions, i.e., scaling functions with the whole support in $I$, is 
\begin{equation}
n_{i}= \left\lfloor \frac{ \left( N_b-1 \right) \, 2^{j(1)-r}-K  }{ b } \right\rfloor + 1,
\end{equation}
where the symbol $\left\lfloor y \right\rfloor$ denotes the largest integer number smaller than $y$.
The number of left boundary scaling functions, i.e., functions that have only a part of the support in the interior of $I$ and their support contains $0$, is 
\begin{equation} \label{number_left_scal}
n_1= K  a-1, \quad a=1/b, 
\end{equation}
and similarly the number of right boundary scaling functions is
\begin{equation} \label{number_right_scal}
n_2= \left\lceil \frac{\left( N_b-1 \right) \, 2^{j(1)-r}}{ b } \right\rceil - \left\lfloor \frac{ \left( N_b-1 \right) \, 2^{j(1)-r}-K }{ b  } \right\rfloor -1.
\end{equation}
Hence, we have
\begin{equation} \label{number_scaling_functions}
\text{ind} \left( 1 \right) = n_1 + n_{i} + n_2 = K  a-1 - \left\lceil \frac{\left( N_b-1 \right) \, 2^{j(1)-r}}{ b } \right\rceil.
\end{equation}
Similarly, the number of wavelet functions on the level $j(l)$ is
\begin{equation} \label{number_wavelets}
\text{ind} \left( 1+l \right) = s a-1 +\left\lceil \frac{ \left( N_b-1 \right) \, 2^{j(l)-r} }{ b } \right\rceil.
\end{equation}
The first $\text{ind}(1)$ columns of $\vDw$ contain values of scaling functions (\ref{scaling_functions_dictionary}), which restricted to $I$ are not identically zero, at points given in (\ref{Nb_points}),
i.e.,
\begin{equation}
    D^W \left( k, l \right) = \phi \left( 2^{j(1)} \frac{k-1}{2^r} - b \left( l - K a  \right)  \right)
\end{equation}
for $k=1, \ldots, N_b$, $l=1, \ldots, \text{ind} (1)$.
The above equation can be recast: 
\begin{eqnarray}
    D^W \left( k, l \right) &=& 
\phi \left( \frac{ \left( k - 1 \right) - b 
    \left( l - K a \right) 2^{r - j(1)} }{ 2^{r - j(1)} } \right) \\
\nonumber    &=& \tphi \left( k - b \left( l - K a \right) 2^{r - j(1)} \right),
\end{eqnarray}
where $\myphi$ is defined by (\ref{def_tphi}) for the level $u=r-j(1)$.
Using the substitution $m=k - b 
    \left( l - K a \right) 2^{r - j(1)}$, we obtain
\begin{equation} \label{dictionary_scaling}
D^W \left( m +  b 
    \left( l - K a \right) 2^{r - j(1)}, l \right) = \tphi(m), \quad m=1, \ldots, K 2^{r-j(1)}+1,
\end{equation}
under the assumption that $1 \leq m +  b 
\left( l - K a \right) 2^{r - j(1)}  \leq N_b$.

The other columns of $\vDw$ contain values of wavelet functions (\ref{wavelets_dictionary}) for levels $j = j(1), \ldots, j(J)$  at points (\ref{Nb_points}), i.e., 
\begin{equation}
    D^W \left( k, n_p + l \right) = \psi \left( 2^{j} \frac{k-1}{2^r} - b \left( l - s a  \right)  \right), \quad s=\frac{M+K}{2},
\end{equation}
for $k=1, \ldots, N_b$, $l=1, \ldots, \text{ind} (j+1)$, and $n_p= \sum_{p=j(1)}^{j} \text{ind} (  p +1 -j(1) )$.
Similarly as above we obtain
\begin{equation} \label{comput_dict_wav}
  D^W \left( m + b \left(l-s a \right) \, 2^{r-j(l)},l+n_p \right)= \tpsi(m),
\end{equation}
where $1 \leq m +  b 
\left( l - s a \right) 2^{r - j(1)} \leq N_b$ and $\mypsi$ is defined by (\ref{def_tpsi}) for the level $u=r-j(1)$.

The following procedure WaveletDict computes a wavelet dictionary.

\begin{center}
{\fontalg 
  \captionsetup{type=figure}
  \captionof{algorithm}{ \\
   Procedure  [$\vDw, \ind, \col$] = WaveletDict(namef, $N_b$, $\j$, $b$)   } \label{algor_WaveletDict}

\begin{algorithmic}

\STATE{{\bf Input:}}\, 

\begin{tabular}{ p{1cm} p{9cm}  }
namef & name of a wavelet family, for available choices see Appendix~A  \\
$N_b$ & number of  points   \\
$\j$ & vector of levels \\
$b$ & translation factor $b=2^{-r_b}$ for some integer $r_b$ 
\end{tabular}

\STATE{{\bf Output:}}\, 

\begin{tabular}{ p{1cm} p{9cm}  }
$\vDw$ & wavelet dictionary  \\
$\ind$ & $ \text{ind}(1)$ is the number of scaling functions at level $j(1)$, and $\text{ind}(k)$ for $k>1$ is the number of wavelets at level $j(k-1)$  \\
$\col$ & cell array such that $\tcol \! \left\{ n \right\} = \left\{ j, k, \text{type}, \text{function} \right\}$ if the $n$-th column of $\vDw$ corresponds to values of a scaling function $\phi(2^j x - b k)$ or a wavelet $\psi(2^j x - b k)$; type=`inner' or `boundary' characterizes type of a function; function=`scaling' or `wavelet' indicates whether the column corresponds to the values of a scaling function or a wavelet\\
\end{tabular}

\vspace{2mm}

\STATE{\COMMENT{Compute scaling and wavelet filters using Algorithm~\ref{algor_Filters} from Appendix~A} }
\STATE{ [$\h$,$\g$,correct\textunderscore name] = Filters(namef) }

\STATE{\COMMENT{Test if a wavelet family namef is available} }

\IF{ correct\textunderscore name$=0$}

\STATE{$\vDw=[ \, ]; \, \ind=[ \, ]; \, \col=[ \, ]; \, \text{return}$ }

\ENDIF
  
\STATE{ $K=\text{length} (\h)-1$}
\COMMENT{support length of $\phi$}

\STATE{ $   s=(K+\text{length} (\g)-1)/2$}
\COMMENT{support length of $\psi$}

\STATE{$r= \left\lceil \log( N_b-1 )/\log(2) \right\rceil$}
\COMMENT{level characterizing $N_b$ (c.f. (\ref{Nb_points})}

\STATE{\COMMENT{Remove levels from $\j$ that contain no inner function} }

\STATE{$j_{min}= \left\lceil \log( s \, 2^r/ (N_b-1) ) / \log(2)) \right\rceil $}
\COMMENT{coarsest possible level}

\STATE{ $\j=\j( \j>= j_{min})$ }
\COMMENT{removing the levels smaller than $j_{min}$}

\COMMENT{Test of parameters}

\STATE{$d_j=\text{length} (\j); \quad r_b= \left\lceil \log(1/b)/ \log(2) \right\rceil$}
\COMMENT{parameter $r_b$ from $b=1/2^{r_b}$}

\IF{$d_j=0$}

\STATE{ fprintf(`no inner functions for these values of levels $\j$, increase $\j$')}

\STATE{$\vDw=[ \, ]; \, \ind=[ \, ]; \, \col=[ \, ]; \, \text{return}$ }

\ELSIF{$r< \max( \j) +r_b$}

\STATE{ fprintf(`small number of points $N_b$ for these values of $\j$ and $b$')}

\STATE{$\vDw=[ \, ]; \, \ind=[ \, ]; \, \col=[ \, ]; \, \text{return}$}
        
\ENDIF

\STATE{\COMMENT{Compute scaling and wavelet generators using Algorithm~\ref{algor_WaveletGen}}}

\STATE{[$\myphi,\mypsi$]=WaveletGen($\h, \g, r-j(1)$)}

\STATE{ \COMMENT{Compute number of scaling functions at level $j(1)$}}

\STATE{$\ind=\text{zeros}(d_j+1,1)$}

\STATE{$\text{ind}(1)=K a-1+ \left\lceil (N_b-1) \, 2^{j(1)-r} /b \right\rceil$}
\COMMENT{c.f. (\ref{number_scaling_functions})}

\STATE{\COMMENT{Compute number of wavelets for level $l$}}

\FOR{$l=1:d_j$}

\STATE{$\text{ind} \left(1+l \right)=s a-1+ \left\lceil (N_b-1) \, 2^{j(l)-r}/b \right\rceil$}
\COMMENT{c.f. (\ref{number_wavelets})}

\ENDFOR

\COMMENT{Compute columns of $\vDw$ corresponding to scaling functions} 

\STATE{$n_f=\text{sum}(\ind); \, \vDw=\text{zeros}(N_b,n_f); \, \col=\text{cell}(n_f,1)$} 

\STATE{$l_s=\text{length}(\myphi); \ n_1=K a-1$ }
\COMMENT{c.f. (\ref{number_left_scal})}

\STATE{$n_2= \left\lceil \left( N_b-1 \right) \, 2^{j(1)-r}  /b \right\rceil - \left\lfloor \left( \left( N_b-1 \right) \, 2^{j(1)-r}-K \right) /b \right\rfloor -1$}
\COMMENT{c.f. (\ref{number_right_scal})}

\STATE{\COMMENT{Compute columns corresponding to inner scaling functions (c.f. (\ref{dictionary_scaling})}}

\FOR{$i=n_1+1:\text{ind}(1)-n_2$}

\STATE{$D^W( b (i-K a) 2^{r-j(1)}+1:b (i-K a) \, 2^{r-j(1)}+K \, 2^{r-j(1)}+1,i)=\myphi$}

\STATE{$\tcol \left\{i \right\}= \left\{ j(1), \,  i-K a, \, \text{`inner'}, \, \text{`scaling'} \right\} $}
        
\ENDFOR

\STATE{\COMMENT{Compute columns corresponding to boundary scaling functions (c.f. (\ref{dictionary_scaling})}}

\FOR{$i=1:n_1$}

\STATE{$D^W(1:l_s-b \, (n_1-i+1) \, 2^{r-j(1)},i)=\text{phi}((n_1-i+1) \, b 2^{r-j(1)}+1:l_s)$ }

\STATE{$ \tcol \left\{i \right\}=\left\{ j(1), -n_1+i-1, \text{`boundary'}, \, \text{`scaling'}
\right\}$ }

\ENDFOR
    
\FOR{ $i=1:n_2$ }

\STATE{$p= \text{ind}(1)-n_2+i$}
\COMMENT{index of column}

\STATE{$D^W(b \, (p-K a) 2^{r-j(1)}+1:N_b,p)=\text{phi}(1:N_b-b(p-K a) 2^{r-j(1)})$ }

\STATE{ $\text{col} \left\{ p \right\}= \left\{ j(1), -n_1+p-1, \text{ `boundary'}, \, \text{`scaling'} \right\} $}
  
\ENDFOR  

\COMMENT{Compute columns of $\vDw$ corresponding to wavelets (c.f. (\ref{comput_dict_wav}))}

\STATE{$n_p=\text{ind}(1)$}
\COMMENT{number of functions on previous levels}

\FOR{$l=1:d_j$}

\STATE{$n_1=s a-1$}

\STATE{$k_1=\left\lfloor ((N_b-1)2^{j(l)-r}-s)/b \right\rfloor, \,
k_2= \left\lceil ((N_b-1)2^{j(l)-r)}/b)-1 \right\rceil $ }

\STATE{ $n_2=k_2-k_1; \,  n_f=n_1+n_2+k_1+1; \, l_w=\text{length}(\mypsi)$}

\FOR{$i=n_1+1:n_f-n_2$}

\STATE{$D^W(b(i-s a)2^{r-j(l)}+1:b (i-s a) 2^{r-j(l)}+s 2^{r-j(l)}+1,i+n_p)=\mypsi$}

\STATE{ $\text{col} \left\{ i+n_p \right\}= \left\{ j(l), i-s a, \text{`inner'}, \, \text{`wavelet'} \right\}$ }

\ENDFOR

\FOR{$i=1:n_1$}

\STATE{$D^W(1: l_w-b (n_1-i+1) 2^{r-j(l)},i+n_p)=\text{psi}( (n_1-i+1) b 2^{r-j(l)}+1:l_w)$}

\STATE{ $\text{col} \left\{ i+n_p \right\}= \left\{ j(l),-n_1+i-1,\text{`boundary'}, \, \text{`wavelet'}  \right\} $}

\ENDFOR

\FOR{$i=1:n_2$}

\STATE{$ p=n_f-n_2+i$}

\STATE{ $D^W (b (p-s a) 2^{r-j(l)}+1:N_b,p+n_p)= \text{psi} (1:N_b-b (p-s a) 2^{r-j(l)})$ }

\STATE{$ \text{col} \left\{ n_p+p \right\}= 
\left\{ j(l), -n_1+p-1, \text{`boundary'}, \, \text{`wavelet'}  \right\} $}

\ENDFOR 

\STATE{ $\mypsi= \text{psi}(1:2:\text{length}(\mypsi)), \, 
n_p=n_p+\text{ind}(l+1)$ }

\ENDFOR

\hrulefill   
\end{algorithmic}
}
\end{center}

The main procedure GenDict validates input parameters,  generates dictionaries $\vDw$ and normalizes their columns. 

\begin{center}
{\fontalg 
  \captionsetup{type=figure}
  \captionof{algorithm}{ 
  \\ Procedure [$\vDw, \ind, \col$]= GenDict(name,pars) } 
\label{algor_GenDict}

\begin{algorithmic}

\STATE{{\bf Input:}}\, 

\begin{tabular}{ p{1cm} p{9cm}  }
name & name of a wavelet family, for available choices see Appendix~A \\
pars & parameters in the form pars = $\left\{  N_b, \j, b  \right\}$  \\
\end{tabular}

\medskip
\STATE{Description of the parameters:

\begin{tabular}{ p{1cm} p{9cm}  }
$N_b$ & number of points  \\
$\j$ & vector of levels  \\
$b$ & translation factor $b=2^{-r_b}$ for some integer $r_b$  \\
\end{tabular}
}

\STATE{{\bf Output:}}\, 

\begin{tabular}{ p{1cm} p{8cm}  }
$\vDw$ & wavelet dictionary  \\
$\ind$ & $ \text{ind}(1)$ is the number of scaling functions at level $j(1)$, and $\text{ind}(k)$ for $k>1$ is the number of wavelets at level $j(k-1)$  \\
$\col$ & cell array such that $\tcol \left\{ n \right\} = \left\{ j, k, \text{type}, \text{function} \right\}$, if the $n$-th column of $\vDw$ corresponds to values of scaling function $\phi(2^j x - b k)$ or wavelet $\psi(2^j x - b k)$; type=`inner' or `boundary' characterizes type of a function; function=`scaling' or `wavelet' indicates whether the column corresponds to the values of a scaling function or a wavelet  \\
\end{tabular}

\vspace{2mm}

\STATE{\COMMENT{Define cell array of names of all available families}}
\STATE{families= $\left\{ \right.$`CW2',`CW3',`CW4',`CDF97',`CDF97d',`CDF53', \\ `Short4',`Short3', `Short2', `Db3',`Db4',`Db5',`Sym3',`Sym4', \\ `Sym5',`Coif26',`Coif38'$\left. \right\}$ }

\STATE{\COMMENT{Validate input parameters} }

\IF{ $\text{nargin} \neq 2$} \STATE{ error(`Need 2 input arguments')} \ENDIF

\IF{ $\sim \! \text{ischar(namef)} $}

\STATE{error(`Name must be a string')}

\ENDIF

\STATE{ $N_b=\text{pars} \left\{1 \right\}; \, \j=\text{pars} \left\{2 \right\}; \, b=\text{pars} \left\{3 \right\}; \j=\text{sort}(\j)$ }

\IF {$b \leq 0$}

\STATE{error(`I expect $b>0$')} 

\ENDIF

\STATE{ $ r= \log(1/b)/ \log(2) $}

\IF{$ \left| r- \text{round} (r) \right| >10^{-10} $}

\STATE{fprintf(`Choose $b$ such that $1/b = 2^r$ for some integer $r$') }

\STATE{$\vDw=[ \, ]; \, \ind=[ \, ]; \, \col=[ \, ]; \, \text{return}$}

\ELSIF{ismember($\left\{ \text{namef} \right\}$,families) }

\STATE{\COMMENT{Generate dictionary using Algorithm~\ref{algor_WaveletDict}}}
\STATE{ [$\vDw, \ind, \col$] = WaveletDict(namef, $N_b$, $\j$, $b$) }

\STATE{\COMMENT{Normalize columns of $\vDw$ using Algorithm~\ref{algor_NormDict} from Appendix~A}}
\STATE{ $\vDw$ = NormDict($\vDw$,1) }
\ELSE

\STATE{  error(`Unknown name of a wavelet family')}

\ENDIF 

\end{algorithmic}
\hrulefill
}
\end{center}

\begin{figure} [!ht] 
 	\includegraphics[height=.23\textheight]{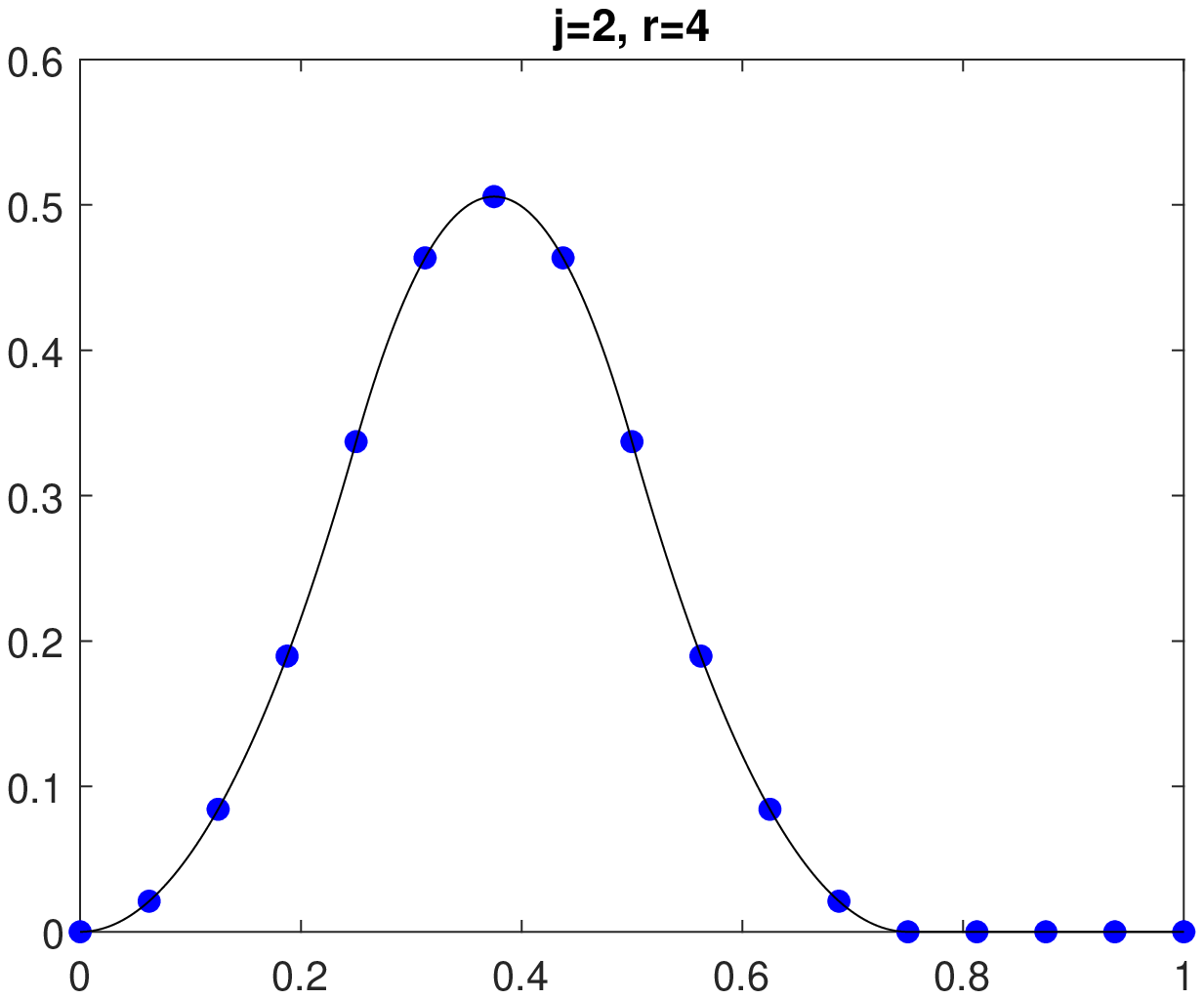}
 	\includegraphics[height=.23\textheight]{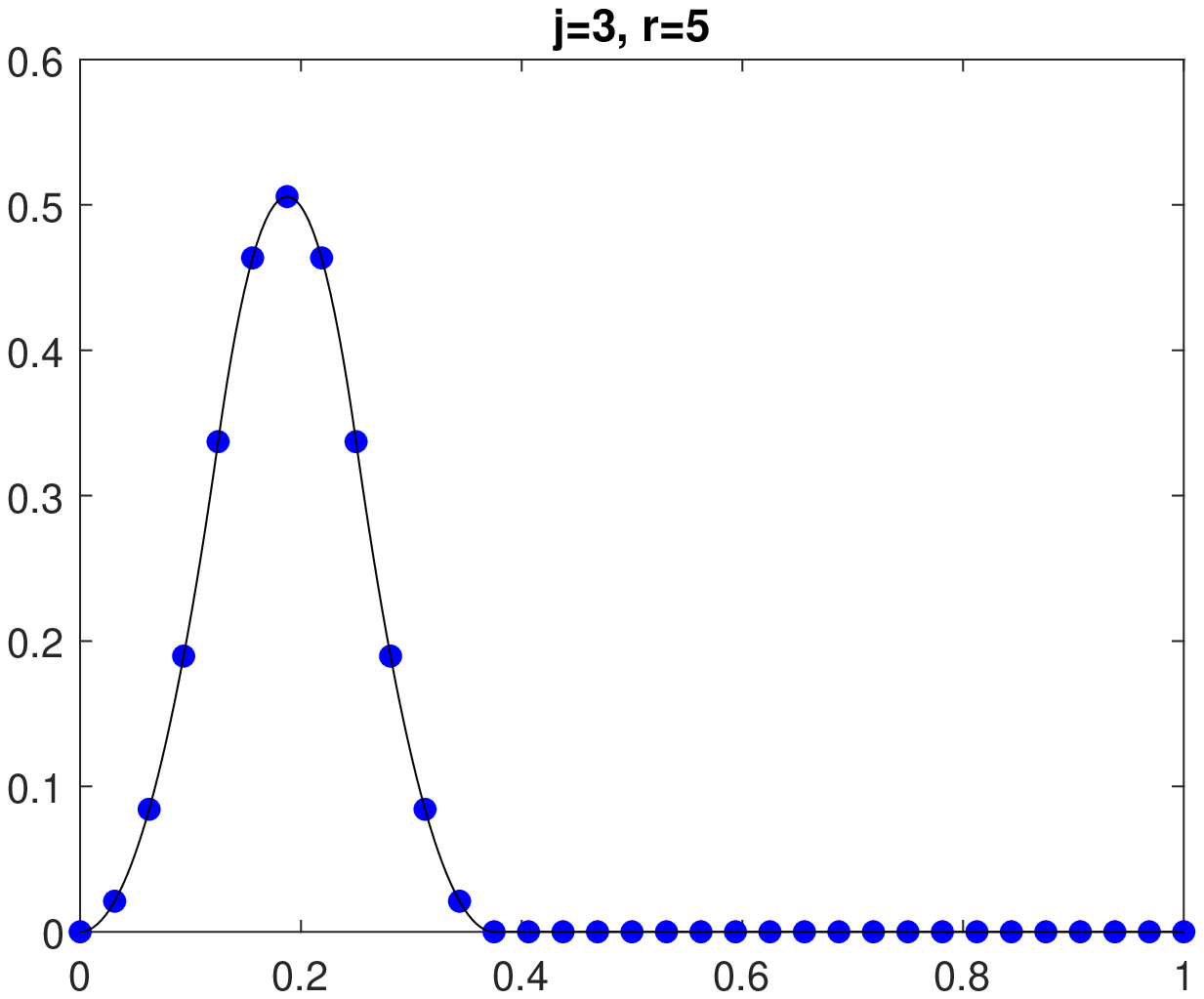}
  \caption{ Vectors of values  $\phi_{j,0,b} \left( l / 2^r \right)$ for $j=2$ and $r=4$ (left) and for $j=3$ and $r=5$ (right).}
  \label{figure2}
\end{figure}

\begin{rem}
It is  worth remarking that the range of scales, say $\j = \left( j_0 , \ldots , J \right)$ depends on length of the signal partition. For a signal segment of length $N_b=2^r+1$ a dictionary contains values of scaling functions and wavelets at points
$l/2^r$ for some integer $r$. For a signal segment of length $2 N_b -1 = 2^{r+1}+1$ a dictionary contains values of functions at points $l/2^{r+1}$. Thus we have
\begin{equation}
 \phi_{j,k,b} \left( \frac{l}{2^r} \right) =   \phi \left( 2^j \frac{l}{2^r} - k b \right) =
    \phi \left( 2^{j+1} \frac{l}{2^{r+1}} - k b \right) = \phi_{j+1,k,b} \left( \frac{l}{2^{r+1}} \right)
\end{equation}
and
\begin{equation}
 \psi_{j,k,b} \left( \frac{l}{2^r} \right) \! = \!   \psi \left( 2^j \frac{l}{2^r} - k b \right) \! = \!
    \psi \left( 2^{j+1} \frac{l}{2^{r+1}} - k b \right) = \psi_{j+1,k,b} \left( \frac{l}{2^{r+1}} \right).
\end{equation}
Therefore, nonzero elements of vectors on the level $j$ in a dictionary for $\R^{N_b}$ correspond to nonzero elements of vectors on the level $j+1$ in a dictionary for $\R^{2 N_b -1}$.
This situation is illustrated in Figure~\ref{figure2}, where vectors of values  $\phi_{j,0,b} \left( l / 2^r \right)$ are displayed for $j=2$ and $r=4$ and for $j=3$ and $r=5$. Note that the nonzero elements in these vectors are the same.
Therefore, if for the signal segment of length $N_b$ the vector $\j = \left( j_0, \ldots, J \right)$ is used, then we recommend to use the vector $\j = \left( j_0, \ldots, J+1 \right)$ for the signal segment of length $2 N_b -1$, and similarly
to use levels $\j = \left( j_0, \ldots, J + m \right)$ for the signal segment of length $2^m N_b -1$.
\end{rem}

\medskip

\begin{figure} [!htb] 
 \includegraphics[height=.3\textheight]{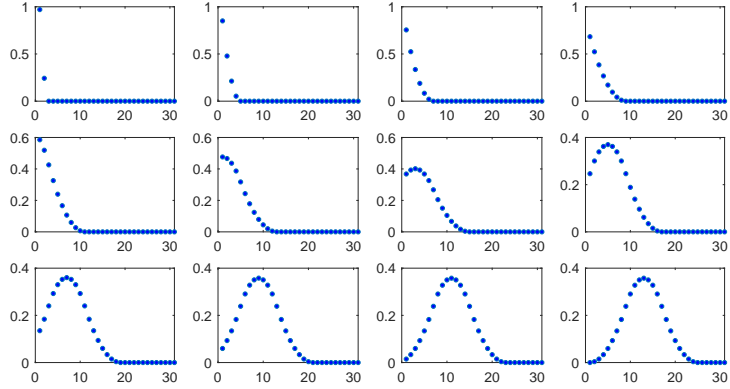}
  \caption{ Plots of $12$ vectors from the dictionary $\vDw$ from Example~1 corresponding to scaling functions on the level $2$. }
  \label{figure3}
\end{figure}

\begin{figure} [!htb] 
 	\includegraphics[height=.3\textheight]{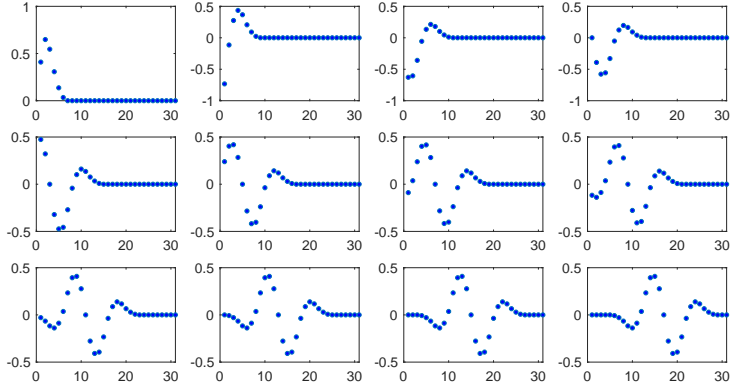}
   \caption{ Plots of $12$ vectors from the dictionary $\vDw$ from Example~1 corresponding to wavelets on the level $2$. }
  \label{figure4}
\end{figure}

\begin{ex}
To build dictionaries for the wavelet family `Short3' at levels 2 and 3, for translation parameter $b = 1/4$, and the number of points $N_b = 33$, use the procedure TestDict below.

\newpage

\begin{center}
\captionsetup[algorithm]{labelformat=empty}
  \captionsetup{type=figure}
\captionof{algorithm}{ Procedure TestDict } 

\begin{algorithmic} 

\vspace{2mm}

\STATE{ namef=`Short3'; $ N_b=33; \, \j=2:3; \, b=1/4 $}

\STATE{ [$\vDw, \ind, \col $]=GenDict(namef,$ \left\{ N_b, \j, b \right\} $)}

\hrulefill
\end{algorithmic}

\end{center}

The output is the matrix $\vDw$ of size $33 \times 97$ and the vector $\ind=[27, 27, 43]$. This means that there are $27$ scaling functions at level 2, 27 wavelets at level 2, and 43 wavelets at level 3. The cell array $\col$ characterizes functions corresponding to columns of $\vDw$. For example 
\begin{equation}
\text{col} \left\{ 30 \right\}=\left\{ 2, -9, \text{`boundary'}, \text{`wavelet'} \right\}
\end{equation}
which means that $30$th column of the matrix $\vDw$ contains values of a wavelet function
$\psi \left( 2^{2} x  - b \left( -9 \right) \right)$.
This wavelet is a boundary wavelet, i.e., only a part of its support lies in the interval $I$ defined by (\ref{def_I}). Some of the vectors from this dictionary corresponding to values of scaling functions are displayed in Figure~\ref{figure3} and some of the vectors corresponding to values of wavelets are displayed in Figure~\ref{figure4}.
\end{ex}

\subsection{Construction of dictionaries for ECG modelling}
\label{sec_dictionaries_ECG}

As mentioned in Sec.~\ref{introduction}, because ECG signals are usually  superimposed to a baseline or smooth background, the full dictionary  $\vD$ we use for  ECG modelling is built  as 
follows
\begin{equation}
\vD=[\vDc \, \vDw],
\end{equation}
where $\vDw$ is the output of Algorithm~\ref{algor_DCos} and $\vDc$ is a matrix containing a few
low frequency components  from a discrete cosine basis. Before normalization $\vDc$ is given as
\begin{equation}
D^C (k,n)  = \cos ( \pi (2k-1) (n-1)/(2 N_b)), \,  k=1, \ldots, N_b, \, n=1, \ldots, M_c, 
\end{equation}
where $M_c$ is a small number in comparison to $N_b$. For
the numerical examples of the next section we consider
$M_c=10$. Algorithm~\ref{algor_DCos} computes $\vDc$. 

\begin{center}
  \captionsetup{type=figure}
  \captionof{algorithm}{  \\
Procedure $\vD^C$ = DCos($N_b, M_c$) } 
\label{algor_DCos}

\begin{algorithmic}

\STATE{ {\bf Input:}}\, 

\begin{tabular}{ p{1cm} p{9.5cm}  }
$N_b$ & the size of the Euclidean space the vectors should belong to \\
$M_c$  & number of frequencies to use starting from 0  \\
\end{tabular}

\STATE{ { \bf Output:}}\, 

\begin{tabular}{ p{1cm} p{9.5cm}  }
$\vD^C$  & matrix whos columns are discrete cosine vectors \\
\end{tabular}

\vspace{2mm}

\STATE{ $n = 1: M_c; \, k = 1:N_b$ }

\STATE{ $\vD^C =  \cos( \pi (2k-1)^T (n-1)/(2 N_b)) $ }

\STATE{ $\vD^C$ = NormDict($\vD^C$,1) }

\end{algorithmic}

\hrulefill

\end{center}

\subsection{Method for construction of the model}
\label{sec_constr_model}

In this section we present the procedures for 
constructing the ECG signal model 
(c.f. Algorithm \ref{ALSM}) and for 
calculating the assessment metrics.
The quality of the signal approximation is assessed
 with respect to the $\PRD$ defined as follows
\be
\label{PRD}
\PRD=\frac{\|\vf - \vfr\|}{\|\vf\|} \times 100 \%, 
\ee
where $\vf$ is the original signal and $\vfr$ is
the signal reconstructed by concatenation 
of the approximated segments $\vfa\cq,\,q=1,\ldots,Q$. 

The local $\PRD$ with respect to every 
 segment in the signal partition is indicated
 as $\prd(q),\,q=1,\ldots, Q$ and calculated as
\be
\label{prd}
\prd(q)=\frac{\|\vf\cq - \vfa\cq\|}{\|\vf\cq\|} \times 100 \%,\quad q=1,\ldots,Q.
\ee
For the signal approximation the 
OOMP method is stopped through a fixed value
 $\tol$ so as to achieve the 
 same value of $\prd$ for all the segments in the
records. Assuming that the target $\prd$ before
quantization is
$\prdo$ we set $\tol = \prdo\|\vf_q\|/100$.

The goal of the signal model is to approximate  each 
segment in the signal partition 
 using as few atoms as possible. Thus, 
for a fixed value of $\PRD$, the sparsity of the signal 
representation is assessed by the sparsity ratio
(SR)
\be
\label{SR}
\text{SR}=\frac{N}{K},
\ee
where $N$ is the total length of the signal
and $K=\sum_{q=1}^Q \kq,$ with $\kq$ the number of 
atoms in the atomic decomposition \eqref{atom} 
of each segment of length $\Nq$.
The corresponding quantity evaluated for 
every cell in the partition is the local sparsity ratio
\be
\label{sr}
\sr(q)=\frac{\Nq}{\kq},\quad q=1,\ldots, Q.
\ee
This local quantity is relevant to the detection of 
non-stationary noise, significant distortion in 
ECG patterns, or changes of morphology in the 
heart beats.

Given an ECG signal $\vf$ the procedure described in
Algorithm $\ref{ALSM}$ constructs the signal
approximation, $\vfr$, using the dictionaries
introduced in the previous section.

\begin{center}
  \captionsetup{type=figure}
  \captionof{algorithm}{  \\
  Procedure [$\vf^{{\rm{r}}},\cl,\vc,\prd,\sr,\PRD,\SR$]=
\text{SignalModel}($\vf,\Nq,\prdo,\name,\pars,\Mc$) } 
  \label{ALSM}

\begin{algorithmic}

\STATE{ {\bf Input:} }\, 

\begin{tabular}{ p{1cm} p{9cm}  }
$\vf$ & signal \\
$\Nq$ & number of points in each segment of the partition \\
$\prdo$ & parameter to control the approximation error\\
$\name$ & name of a wavelet family \\
$\pars$ & parameters as described in Algorithm~\ref{algor_GenDict}\\
$\Mc$ & number $\Mc$ of components in the cosine subdictionary \\
\end{tabular}

\STATE{ {\bf Output:} }\, 

\begin{tabular}{ p{1cm} p{9cm}  }
$\vf^{{\rm{r}}}$ & approximated signal \\
$\cl$ & cell with the indices of the atoms in the 
atomic decomposition of each element in the partition \\
$\vc$ & cell with the coefficients in the atomic decomposition of each element in the 
partition \\
$\prd$ & vector $\prd \in \R^Q$ (cf. \eqref{prd})  \\
$\sr$ & vector $\sr \in  \R^Q$ (cf. \eqref{sr})\\
$\PRD$ & global PRD \\
$\SR$ & global SR \\
\end{tabular}

\STATE{\COMMENT{Create the signal partition using Algorithm~\ref{ALSP}}}

\STATE{[$\vfc,Q,\vf$]=\text{Partition}($\vf,\Nq$)}

\STATE{\COMMENT{Construct the wavelet dictionary $\vDw$ 
using Algorithm~\ref{algor_GenDict} given in Appendix~B}}

\STATE{[$\vDw$, ind]= GenDict($\name,\pars$)}

\STATE{\COMMENT{Construct the component $\vDc$  
 using  Algorithm~\ref{algor_DCos} given in Appendix~B}}

\STATE{[$\vDc$]=DCos($\Nq,\Mc$)}

\STATE{\COMMENT{Merge $\vDc$ and $\vDw$ to create 
dictionary $\vD$}}

\STATE{$\vD= [\vDc \,\, \vDw]$}
 
\STATE{Set $\vfr=[\;]$,  $K=0$ and $N=$length$(\vf)$.}
\FOR {q=1:Q}
\STATE{$\tol=\prdo\|\vfc\cq\|/100$}
\STATE{\COMMENT{Call the OOMP function to construct the model (c.f. \eqref{atom})}}
\STATE{[$\vfa\cq,\cl\cq,\vc\cq$]=OOMP($\vfc\cq,\vD,\tol,1$)}
\STATE{\COMMENT{Calculate local sr and prd (c.f. \eqref{sr} and  \eqref{prd}}}
\STATE{$\prd(q)=\frac{\|\vfc\cq-\vfa\cq\|}{\|\vfc\cq\|}\times 100$}
\STATE{$k(q)=\text{length}(\vc\cq))$}
\STATE{$\sr(q)=\Nq/k(q)$}
\STATE{$K=K+k(q)$}
\STATE{$\vfr=[\vfr \,\vfa\cq]$}
\ENDFOR
\STATE{\COMMENT{Calculate global SR and PRD (c.f. \eqref{SR} 
and \eqref{PRD})}}
\STATE{$\SR=N/K$; $\PRD=\frac{\|\vf -\vfr\|}{\|\vf\|} \times 100$}
\end{algorithmic}
\hrulefill
\end{center}

\section{Results }
We illustrate now the use of the software 
to approximate records 117, 202, and 231 
in the MIT-BIH Arrhythmia database. 
Each record consists of   
650000 samples and is partitioned for the 
approximation in segments of $\Nq=500$ points 
each. Table~\ref{Table1} gives the values of 
the SR (c.f.~\eqref{SR}) 
achieved using wavelet bases, denoted as 
$\SRB$, and  wavelet dictionaries  
 denoted as $\SRD$. The wavelet families are
indicated in the first column of Table~\ref{Table1}. 
The wavelet dictionary is constructed with 
scales $\j= \left( 3,\ldots,7 \right)$ and translation parameter $b=1/4$, 
whilst the wavelet basis entails to add one more 
scale and a translation parameter 
$b=1$. In all the cases the approximation 
is realized to obtain $\PRD=0.51\%$. 

Table \ref{Table1} is produced by running the script
 `Run\_ECG\_Appox' and changing the variable `namef' 
 to the corresponding family option. 
 
 \begin{center}
\captionsetup[algorithm]{labelformat=empty}
  \captionsetup{type=figure}
\captionof{algorithm}{ Procedure Run\textunderscore ECG\textunderscore Approx  } 
\label{AlExam}
\begin{algorithmic}
\STATE{\COMMENT{Read the signal $\vf$}}
\STATE{file=`Record\_231\_11bits.dat'}
\STATE{fid=fopen(file,`r')}
\STATE{$\vf$=fread(fid,`ubit11')}
\STATE{fclose(fid)}
\STATE{\COMMENT{Set the required PRD for the approximation}}
\STATE{$\prdo=0.53$}
\STATE{\COMMENT{Set the length for each segment in the signal partition}}
\STATE{$\Nq=500$}
\STATE{\COMMENT{Set the parameters for the wavelet dictionary}}
\STATE{$\name=$`CDF97'; $ b=0.25; \, \j = 3:7; \, \pars=\{ N_b, \j, b\}$}
\STATE{\COMMENT{Set the number of cosine components}}
\STATE{$\Mc=10$}
\STATE{\COMMENT{Construct the signal module}}
\STATE{[$\vf^{{\rm{r}}},\cl,\vc,\prd,\sr,\PRD,\SR$]= \text{SignalModel}($\vf,\Nq,\prdo,\name,\pars,\Mc$)}
\STATE{\COMMENT{Plote the first 2000 sample points in the signal, the approximation and the error}}
\end{algorithmic}
\hrulefill
\end{center}
 
\def\CWC{\mathbf{CW4}}
\def\CWQ{\mathbf{CW3}}
\def\CWL{\mathbf{CW2}}
\def\CDF{\mathbf{CDF97}}
\def\CDFd{\mathbf{CDF97d}}
\def\CDFF{\mathbf{CDF53}}
\def\Dbq{\mathbf{Db3}}
\def\Dbc{\mathbf{Db4}}
\def\Dbf{\mathbf{Db5}}
\def\Symq{\mathbf{Sym3}}
\def\Symc{\mathbf{Sym4}}
\def\Symf{\mathbf{Sym5}}
\def\Shl{\mathbf{Short2}}
\def\Shq{\mathbf{Short3}}
\def\Shc{\mathbf{Short4}}
\def\Coifts{\mathbf{Coif26}}
\def\Coifte{\mathbf{Coif38}}

\begin{table} [htb!]
\caption{SRs achieved using dictionaries with 
$\vDw$ component as indicated in the first column 
 of the table. $\SRB$ are values of SR obtained if 
$\vDw$ is a basis and $\SRD$ if $\vDw$ is a dictionary.}
\label{Table1}
\begin{center}
\begin{tabular}{|l||r|r|| r| r|| r| r||}
\hline
Rec. & \multicolumn{2}{|c||}{117} & 
\multicolumn{2}{|c||}{202} & \multicolumn{2}{|c||}{231}\\ \hline \hline
$\vDw$ & $\SRB$ & $\SRD$& $\SRB$ & $\SRD$& $\SRB$& $\SRD$\\ 
\hline  \hline
$\CWL$&17.5& 26.5& 17.3&24.5& 15.7 & 23.0 \\ 
$\CWQ$& 17.4 & 28.1& 15.9& 24.9 & 15.6 & 24.0   \\
$\CWC$&15.7& 24.8&  14.3& 22.5&18.4&21.9 \\
$\CDF$&21.5& 30.3&  21.4& 28.4&19.5&27.5 \\
$\CDFd$&17.2 & 23.5 & 17.3 & 22.5& 15.8& 21.9\\
$\CDFF$&22.4 & 29.6& 23.6& 27.0& 20.2& 27.0 \\
$\Dbq$&18.5 & 23.7& 18.1& 22.7& 16.6& 22.7\\
$\Dbc$&19.0 & 25.7& 19.1&24.7 & 17.7& 24.1 \\
$\Dbf$&20.4 & 26.1& 18.7&24.2 & 17.8 & 24.1 \\
$\Symq$&18.4 & 23.8& 18.1&22.7 & 16.6 & 22.7\\
$\Symc$&19.7 & 27.5& 19.5& 25.8 & 17.7 & 25.1 \\
$\Symf$&20.5& 28.3& 20.6& 28.5 & 18.4 & 25.4 \\
$\Shl$&8.2& 27.9&  8.7& 26.3 &  8.1 & 24.7 \\
$\Shq$&19.6&31.8& 18.3& 27.6 &  17.8&27.3 \\
$\Shc$&9.5 & 29.1& 10.1 & 27.6 & 9.1 & 26.6 \\
$\Coifts$&17.7& 23.0& 17.7 & 21.8 & 16.3 &24.7 \\ 
$\Coifte$&19.5& 28.5 &19.7 & 26.5 &17.8 & 26.1  \\ 
\hline \hline
\end{tabular}
\end{center}
\end{table}



\begin{figure}[htb!]
\begin{center}
\includegraphics[width=6cm]{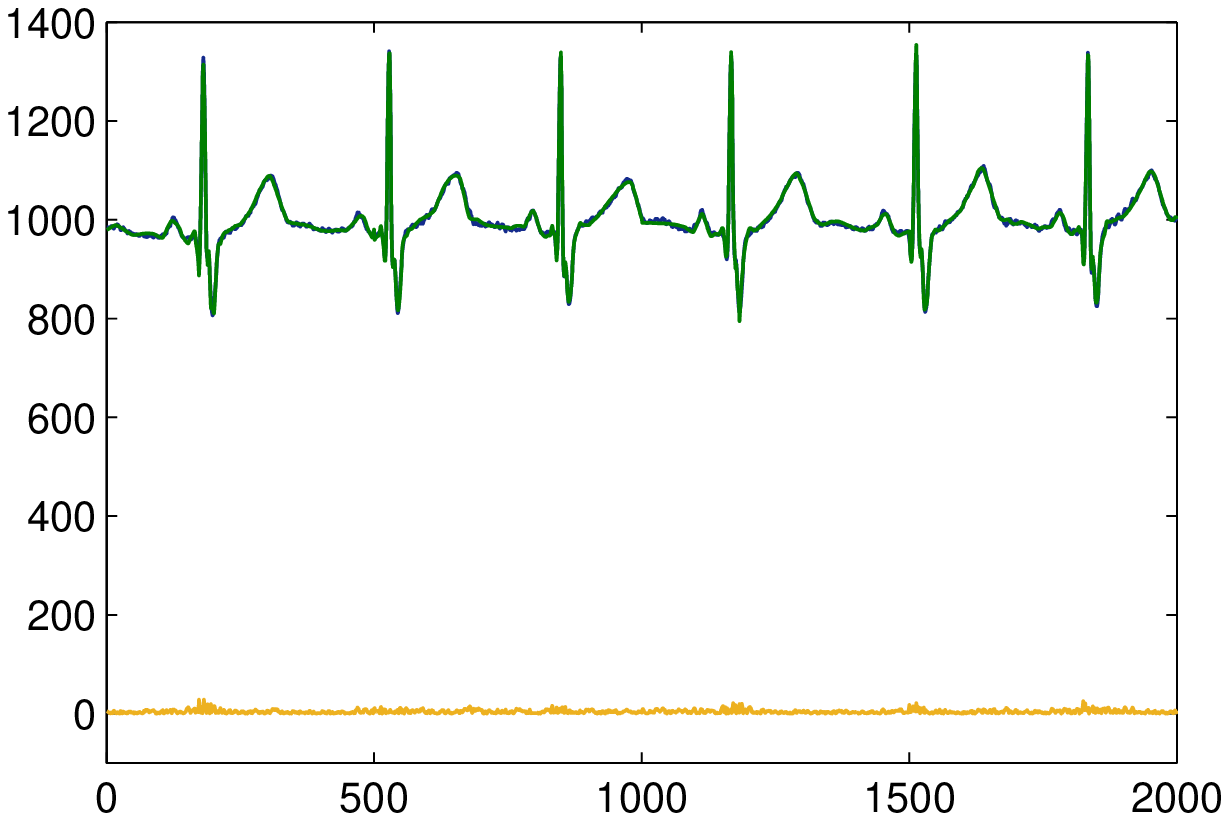}
\includegraphics[width=6cm]{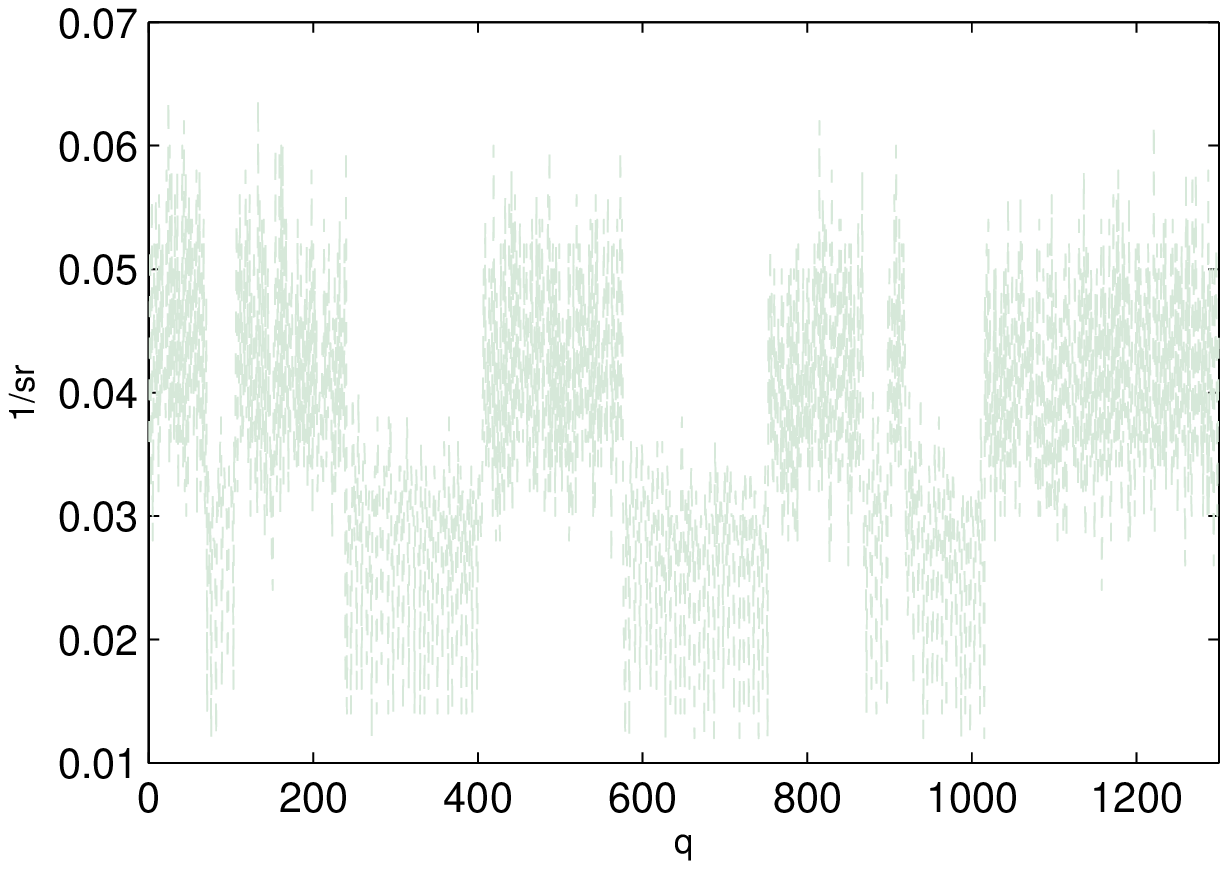}\\
\includegraphics[width=6cm]{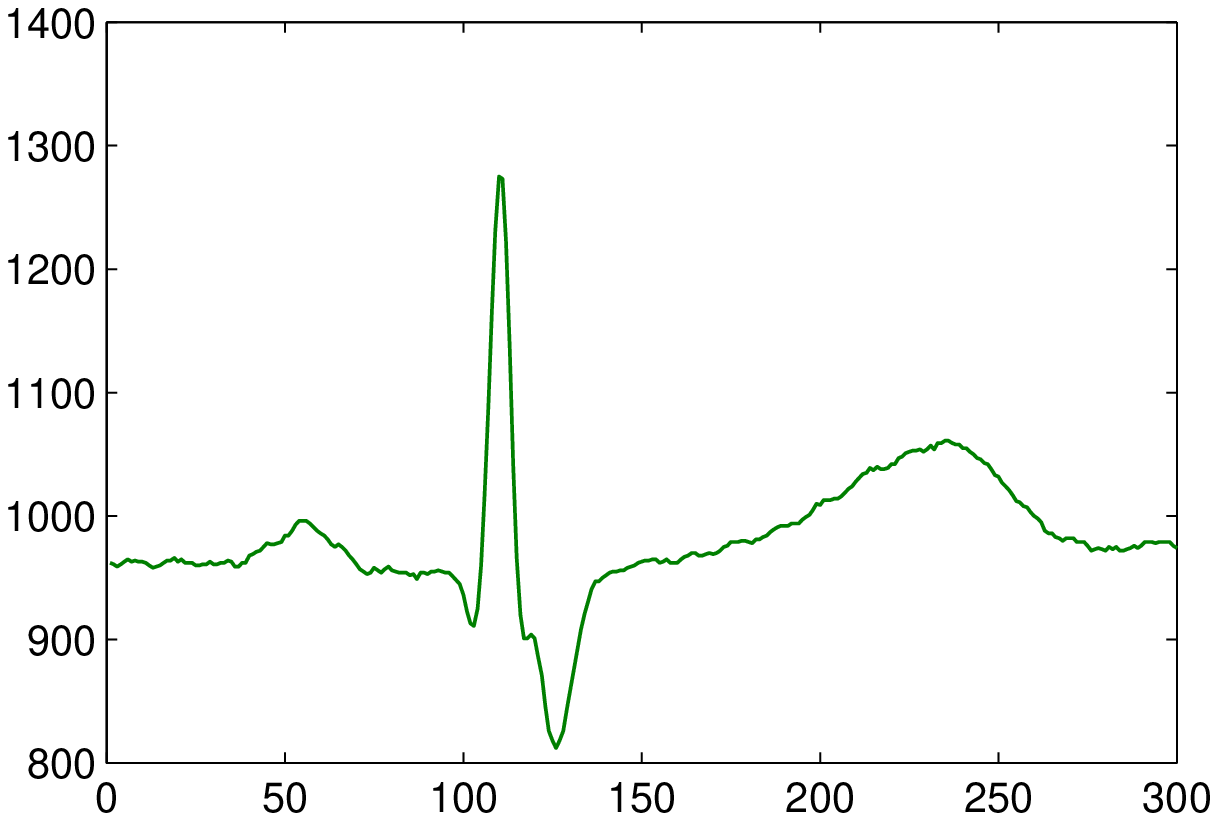}
\includegraphics[width=6cm]{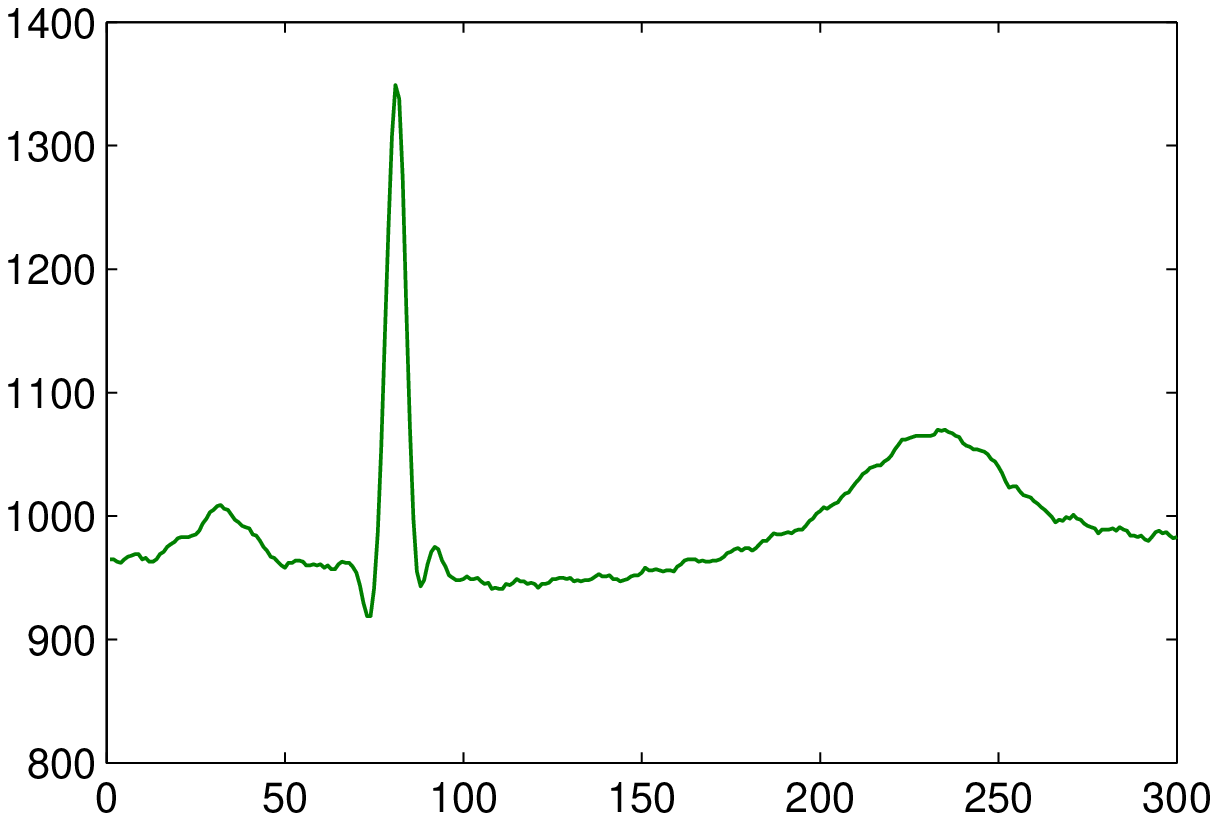}
\caption{The waveforms in the top left graph
 are the raw data and the approximations corresponding
to 2000 points in the records 231 (the bottom 
line in the same graph is the point-wise error).  The 
top right graph plots the
values $1/\sr$ for records 231. 
The bottom left graph is a typical heart beat in a segment 
for which the values of $1/\sr(q)$  belongs to 
 the upper band. In the bottom  
right graph the heart beat corresponds to a frame in the 
lower band.  
}
\label{figure5}
\end{center}
\end{figure}

The top left graph in Figure~\ref{figure5} illustrates the first 2000 
 points in the record 231 and the approximation for 
$\PRD=0.51\%$. The top right graph represents the 
values of local sparsity $1/\sr(q)$, $q=1,\ldots,1300$  
for the same record. It is noticed that these 
values can be classified into two well defined bands. The 
 bottom left graph in Figure~\ref{figure5} shows a typical 
heart beat in a frame  corresponding to a 
value $1/\sr$ in the upper band, and the 
 bottom right graph to a value in the lower band.
The  morphologic difference between the two heart beats 
is noticeable at a glance. 

\section{Discussion}

As observed in the Table \ref{Table1}, the gain in dimensionality reduction (larger value of SR) is significant when considering a wavelet dictionary, instead of a wavelet basis, as component $\vDw$ of the full dictionary $\vD$. This result was demonstrated in \cite{RNC19} on the whole MIT-BIH Arrhythmia database, which motivated the present work to provide the details and algorithms for the actual construction of dictionaries from different wavelet prototypes. Notice that dictionaries for the families CDF97, CDF53, and Short3 produce the highest sparsity ratios. This is also in line with the results presented in \cite{RNC19} for the  MIT-BIH Arrhythmia data set.

In \cite{RNC19} wavelet dictionaries have been shown to be suitable for lossy compression at low level distortion.  However, dimensionality reduction is also useful for other applications. It is envisaged that the model
could be relevant to analysis and classification tasks \cite{RR18}.

While the results have been obtained using the  OOMP  approach for selecting the elementary components 
in the model,  other selection strategies \cite{DTD06, EKB10, NT09, LRN16} could be applied with the identical dictionaries. The focus of this work was the  construction of dictionaries delivering piecewise sparse approximation of ECG signals using any
suitable approach for the selection process.

\section{Conclusions}

A detailed description of methods, algorithms, and usage of the software for the construction of wavelet dictionaries has been presented. The use of the software, which has been made publicly available on a dedicated website \cite{RNC19a}, was illustrated to reduce the dimensionality of three records from the MIT-BIH Arrhythmia database. For all the wavelet families, the sparsity ratio yielded by dictionaries with translation parameter $b = 1/4$  was shown to be significantly larger than for the corresponding wavelet bases. The conclusions
coincide with those that were drawn in the previous publication \cite{RNC19} using the whole database. The purpose of this paper was to provide a complete description of the construction of the wavelets dictionaries, which had not been addressed in \cite{RNC19}.  We believe the proposed dictionaries should be of assistance to general applications which relay on dimensionality reduction, at low level distortion, as a first step of further ECG signal processing.

\section*{Declaration of Competing Interest}

There are no known conflicts of interest associated with this publication. 

\section*{Statement of Ethical Approval}

Ethical approval is not required for this work.

\section*{Acknowledgements}

This research did not receive any specific grant from funding agencies in the public, commercial, or not-for-profit sectors. 

\section*{Appendix A}

In this appendix, we present auxiliary procedures used in algorithms in Section~\ref{sec_wavelet_dict}. In the algorithms, `namef' denotes a name of a wavelet family, available choices are:
 
\medskip

\begin{longtable}{ r @{=}p{2.3cm}  p{7.6cm} }
 namef & `CW2' & Chui-Wang linear spline wavelets \cite{CW92} \\
       & `CW3' & Chui-Wang quadratic spline wavelets \cite{CW92} \\
       & `CW4' & Chui-Wang cubic spline wavelets \cite{CW92} \\
       & `CDF97' & primal CDF97 wavelets \cite{CDF92} \\
       & `CDF97d'   & dual CDF97 wavelets \cite{CDF92} \\
       & `CDF53' & primal CDF53 wavelets \cite{CDF92} \\
       & `Short4'       & cubic spline wavelet with short support and 4 vanishing moments \cite{Chen95, Han06} \\
       & `Short3'       & quadratic spline wavelet with short support and 3 vanishing moments \cite{Chen95, Han06} \\
       & `Short2'       & linear spline wavelet with short support and 2 vanishing moments \cite{Chen95, Han06} \\
       & `Db3'          & Daubechies wavelet with 3 vanishing moments \cite{Daub88} \\
       & `Db4'          & Daubechies wavelet with 4 vanishing moments \cite{Daub88} \\
       & `Db5'          & Daubechies wavelet with 5 vanishing moments \cite{Daub88} \\
       & `Sym3'         & symlet with 3 vanishing moments \cite{Daub93}  \\
       & `Sym4'         & symlet with 4 vanishing moments \cite{Daub93}  \\
       & `Sym5'         & symlet with 5 vanishing moments \cite{Daub93}  \\
       & `Coif26'       & coiflet with 2 vanishing moments and the support length 6 that is most regular \cite{Daub93} \\ 
       & `Coif38'       & coiflet with 3 vanishing moments and the support length 8 that is most symmetrical  \cite{Daub93} \\ 
\end{longtable}

A wavelet basis is determined by its scaling and wavelet filters. Algorithm~\ref{algor_Filters} assigns these filters for a chosen wavelet family, the values of filters are computed by methods from \cite{ CDF92, Chen95, CW92, CFK08, Daub88, Daub93, Han06}.

\begin{center}
  \captionsetup{type=figure}
  \captionof{algorithm}{ \\  Procedure [$\h$,$\g$,correct\textunderscore name] = \text{Filters}(namef) } \label{algor_Filters}

\begin{algorithmic} 

\STATE{ {\bf Input:} } \, 

\begin{tabular}{ p{2cm} p{8cm}  }
namef & name of a wavelet family \\
\end{tabular}

\STATE{ {\bf Output:} } \, 

\begin{tabular}{ p{2cm} p{8cm}  }
$\h$ & scaling filter for a wavelet family specified by `namef'  \\
$\g$ & wavelet filter for a wavelet family specified by `namef'  \\
correct\textunderscore name & returns $1$ if `namef' is a name of an available wavelet family, otherwise returns $0$ \\
\end{tabular}

\vspace{2mm}

\STATE{ correct\textunderscore name=1 }

\SWITCH{namef}

\CASE {`CW2'}
  \STATE{ $\h=[1/2, 1, 1/2]; \,  \g=[1, -6, 10, -6, 1]/12$}
\ENDCASE

\CASE{ `CW3' } 
   \STATE{ $\h=[1/4, 3/4, 3/4, 1/4]; \, \g=[1, -29, 147, -303, 303, -147, 29, -1]/480$ }
\ENDCASE

\CASE{ `CW4' }
\STATE{ $\h=[1/8, 1/2, 3/4, 1/2, 1/8]$ }
\STATE{ $\g=[1, -124, 1677, -7904, 18482, -24264, 18482, -7904, 1677, -124, 1]/2520$ }
\ENDCASE        
    
\CASE{`CDF97'}
\STATE{ $\h=[-0.045635881557,  -0.028771763114,  0.295635881557,  0.557543526229,$ \\ $0.295635881557, -0.028771763114, -0.045635881557]$}

\STATE{ $\g=[ 0.026748757411, 0.016864118443, -0.078223266529, -0.266864118443$ \\
$0.602949018236, -0.266864118443, -0.078223266529, 0.016864118443,$ \\ $0.026748757411]$ }
\ENDCASE    

\CASE{`CDF97d'}
\STATE{$ \h=[ 0.026748757411, -0.016864118443, -0.078223266529, 0.266864118443,$\\
$ 0.602949018236, 0.266864118443, -0.078223266529, -0.016864118443,$ \\
$0.026748757411000]$ }
\STATE{$ \g=[0.045635881557, -0.028771763114, -0.295635881557,  0.557543526229,$ \\ $-0.295635881557, -0.028771763114, 0.045635881557]$ }
\ENDCASE

\CASE{`CDF53'}
\STATE{$\h=[1/2, 1, 1/2]; \, \g=[-1/8, -1/4, 3/4, -1/4, -1/8]$}
\ENDCASE

\CASE{ `Short4'}  
\STATE{$\h=[1/8, 1/2, 3/4, 1/2, 1/8]; \, \g=[1/8, -1/2, 3/4, -1/2, 1/8]$}
\ENDCASE

\CASE{ `Short3'}
\STATE{$ \h=[1/4, 3/4, 3/4, 1/4]; \, \g=[-1/4, 3/4, -3/4, 1/4]$}
\ENDCASE

\CASE{`Short2'}
\STATE{ $\h=[1/2, 1, 1/2]; \, \g=[-1/2, 1, -1/2]$ }
\ENDCASE

\CASE{ `Db3' }
\STATE{ $\h = [0.035226291882101,  -0.085441273882241,  -0.135011020010391,$ \\   $ 0.459877502119331,   0.806891509313339,   0.332670552950957]$}
\STATE{$\g = [-0.332670552950957,   0.806891509313339,  -0.459877502119331,$ \\   $-0.135011020010391,   0.085441273882241,   0.035226291882101]$}
\ENDCASE

\CASE{`Db4'}
\STATE{$ \h=[0.162901714025620, 0.505472857545650, 0.446100069123190,$ \\ $-0.019787513117910, -0.132253583684370, 0.021808150237390,$ \\
$0.023251800535560, -0.007493494665130]$ }
\STATE{$ \g=-\text{fliplr}([0.162901714025620, -0.505472857545650, 0.446100069123190,$ \\
$0.019787513117910,  -0.132253583684370, -0.021808150237390,$ \\
$0.023251800535560, 0.007493494665130])$ }
\ENDCASE

\CASE{`Db5'}
\STATE{ $\h=[ 0.003335725285002, -0.012580751999016, -0.006241490213012,$ \\
$0.077571493840065, -0.032244869585030,  -0.242294887066190,$ \\   $0.138428145901103, 0.724308528438574,  0.603829269797473,$ \\   $0.160102397974125]$ }
\STATE{ $\g=[ -0.160102397974125, 0.603829269797473, -0.724308528438574,$ \\
$0.138428145901103, 0.242294887066190, -0.032244869585030,$ \\
$-0.077571493840065, -0.006241490213012, 0.012580751999016,$ \\
$0.003335725285002]$ }
\ENDCASE

\CASE{`Sym3'}
\STATE{ $\h=[0.035226291882101,  -0.085441273882241,  -0.135011020010391,$ \\
$ 0.459877502119331,   0.806891509313339,   0.332670552950957]$ }
\STATE{$\g=[-0.332670552950957,   0.806891509313339,  -0.459877502119331,$ \\
  $ -0.135011020010391, 0.085441273882241, 0.035226291882101]$ }
\ENDCASE
    
\CASE{ `Sym4'}
\STATE{$\h=[0.022785172948000,  -0.008912350720850, -0.070158812089500,$ \\   $0.210617267102000, 0.568329121705000, 0.351869534328000,$ \\
$-0.020955482562550,  -0.053574450709000]$ }
\STATE{$\g=\text{fliplr}([0.022785172948000,  0.008912350720850, -0.070158812089500,$ \\
$-0.210617267102000, 0.568329121705000,  -0.351869534328000,$ \\
$-0.020955482562550,  0.053574450709000])$ }
\ENDCASE

\CASE{ `Sym5' }
\STATE{$\h=[0.027333068345078,   0.029519490925775,  -0.039134249302383,$\\    $0.199397533977394, 0.723407690402421,   0.633978963458212,$ \\  
$0.016602105764522,  -0.175328089908450,  -0.021101834024759,$ \\  $0.019538882735287]$ }
\STATE{ $\g=[ -0.019538882735287,  -0.021101834024759,   0.175328089908450,$ \\
$0.016602105764522, -0.633978963458212,  0.723407690402421,$ \\  $-0.199397533977394,  -0.039134249302383,  -0.029519490925775,$ \\   $0.027333068345078]$ }
\ENDCASE

\CASE{`Coif26'}
\STATE{$\h=[9- \sqrt(15), 13+ \sqrt(15), 6+2 \sqrt(15), 6-2 \sqrt(15), 1- \sqrt(15),$ \\
$-3+ \sqrt(15)]/32$ }
\STATE{$\g=-\text{fliplr}([9- \sqrt(15), -13- \sqrt(15), 6+2 \sqrt(15), -6+2 \sqrt(15),$ \\
$1- \sqrt(15), 3- \sqrt(15)]/32)$ }
\ENDCASE

\CASE{`Coif38'}
\STATE{$\h=[-1/32- \sqrt(7)/128, -3/128, 9/32+3 \sqrt(7)/128, 73/128, 9/32-3 \sqrt(7)/128, -9/128, -1/32+ \sqrt(7)/128, 3/128]$ }
\STATE{$\g=-\text{fliplr}([-1/32- \sqrt(7)/128, 3/128, 9/32+3 \sqrt(7)/128, -73/128, 9/32-3 \sqrt(7)/128, 9/128, -1/32+ \sqrt(7)/128, -3/128])$ }
\ENDCASE

\OTHERWISE

\STATE{ disp(`wrong name of a wavelet family') }

\STATE{ correct\textunderscore name = 0 }

\ENDOTHERWISE
        
\ENDSWITCH    
  
\hrulefill
\end{algorithmic}

\end{center}

Now, we introduce a simple procedure NormDict for normalization of dictionaries. More precisely, this procedure normalizes the columns of dictionary $\vD$ to have the Euclidean norm equaled to $1 / \sqrt{ \delta }$.

\begin{center}
  \captionsetup{type=figure}
  \captionof{algorithm}{  \\ Procedure 
  $\vD$ = NormDict($\vD$, $\delta$) } \label{algor_NormDict}

\begin{algorithmic}

\STATE{ {\bf Input:} } \, 

\begin{tabular}{ p{1cm} p{9cm}  }
$\vD$ & wavelet dictionary \\
$\delta$ & parameter such that prescribed norm size is $1/ \sqrt{ \delta }$
\end{tabular}

\STATE{ {\bf Output:} } \, 

\begin{tabular}{ p{1cm} p{9cm}  }
$\vD$ & normalized wavelet dictionary such that the Euclidean norm of each column is $1/ \sqrt{ \delta }$ \\
\end{tabular}

\vspace{2mm}

\STATE{tol=$10^{-5}$ }  

\IF {nargin=1} 

\STATE{  $\delta=1$ }

\ENDIF   

\STATE{$N$=size($\vD$,2); $i=0$}

\WHILE{ $ i <N$ }

\STATE{ $i=i+1; \,  \text{nor}= \sqrt{ \delta } \, \left\| D(:,i) \right\| $}

\IF { $ \text{nor} > \text{tol}$ } 

\STATE{ $D(:,i)= D(:,i)/ \text{nor}$ }

\ELSE  

\STATE{ $D(:,i)=[ \, ];  N=N-1$ }

\ENDIF

\ENDWHILE

\hrulefill
\end{algorithmic}

\end{center}

\section*{Appendix B}

In this appendix, we present auxiliary procedures used in algorithms in Section~\ref{sec_constr_model}. 
The next procedure Partition creates a partition of the signal
$\vf$ into $Q$ segments of the prescribed length $\Nq$.

\begin{center}
  \captionsetup{type=figure}
  \captionof{algorithm}{ \\
Procedure [$\vfc,Q,\vf$]=\text{Partition}($\vf,\Nq$) } 
\label{ALSP}
 
\begin{algorithmic}

\STATE{ {\bf Input:} }\,

\begin{tabular}{ p{1cm} p{9cm}  }
$\vf$ & signal \\
$\Nq$ & length of each segment in the partition \\
\end{tabular}

\STATE{ {\bf Output:} }\,

\begin{tabular}{ p{1cm} p{9cm}  }
$\vfc$ & cells $\vfc\cq,\,q=1,\ldots,Q$ with the signal partition \\
$Q$ & number of cells in the partition  \\
$\vf$ & resized signal to be of length $Q\Nq$ \\
\end{tabular}

\STATE{$N=\text{length($\vf$)}; Q=\lfloor \frac{N}{\Nq} \rfloor; t_o=1$}

\STATE{$\vf \leftarrow f(1: Q \Nq)$}

\FOR {$q=1:Q$}
\STATE{$t=t_o:t_o +\Nq-1; \, t_0=t_0 + \Nq$}
\STATE{$\vfc\cq=f(t)$}
\ENDFOR
\end{algorithmic}
\hrulefill
\end{center}

The procedure for signal approximation using OOMP method is presented below. 

\begin{center}
  \captionsetup{type=figure}
\captionof{algorithm}{ \\ Procedure [$\vf^{{\rm{a}}},\cl,\vc$]=
\text{OOMP}($\vf,\vD,\tol,\ell_1$) } 
\label{ALOOMP}

\begin{algorithmic}
\STATE{ {\bf Input:} }\, 

\begin{tabular}{ p{1cm} p{9cm}  }
$\vf$ & signal to be approximated by an atomic decomposition \\
$\vD$ & wavelet dictionary  \\
$\tol$ & parameter to control the approximation error\\
$\ell_1$ & index of the atom for initializing the OOMP algorithm \\
\end{tabular}

\STATE{ {\bf Output:} } \, 

\begin{tabular}{ p{1cm} p{9cm}  }
$\vf^{{\rm{a}}}$ &  approximation of the signal $\vf$ (c.f. \eqref{atom}) \\
$\cl$ & vector whose components are the indices of the selected columns from the input dictionary  \\
$\vc$  & coefficients $\vc \in \R^{\Nq}$ of the atomic decomposition
(c.f. \eqref{atom}) \\
\end{tabular}

\STATE{\COMMENT{The method implemented in this function
is fully described in the main paper \cite{RNC19}.}}

\end{algorithmic}
\hrulefill
\end{center}


\end{document}